\documentclass{amsart}

\usepackage[T1]{fontenc}
\usepackage{graphicx} 
\usepackage{amsmath,amsfonts,amssymb,amsthm}
\usepackage{tikz}
\usetikzlibrary{arrows,automata}
\usepackage{tikz-cd}
\usepackage{hyperref}
\usepackage{placeins}
\usepackage{bm}
\usepackage{thmtools}
\usepackage{mathtools}
\usepackage{mathptmx}
\usepackage{tikz,hyperref,stmaryrd,a4wide,amssymb,enumerate}

\usepackage[utf8]{inputenc}     
\usepackage{mathtools}
\usepackage{amsfonts}
\usepackage{amsthm}
\usepackage{amssymb}
\usepackage{amsmath}
\usepackage{tikz}
\usetikzlibrary{arrows,automata}
\usepackage{stmaryrd}
\usepackage{cite}
\usepackage{amsmath}
\usepackage{hyperref}
\usepackage{mathtools}
\usepackage{xcolor}

\calclayout

\usepackage{soul} 

\newcommand{\N}{\mathbb{N}}

\newcommand{\Z}{\mathbb{Z}}
\newcommand{\R}{\mathbb{R}}
\newcommand{\C}{\mathbb{C}}

\newcommand{\ii}{\mathrm{i}}
\DeclareMathOperator{\lcm}{lcm}

\def\claimmargin{1.5em}

\newenvironment{claim}{\list{} {\rightmargin\claimmargin \leftmargin\claimmargin} \item[] \textbf{Claim.}\hspace{0.5em}\begin{itshape}}{\end{itshape}\endlist}

\newenvironment{claimproof}{\list{} {\rightmargin\claimmargin \leftmargin\claimmargin} \item[]\textit{Proof.}\hspace{1em}}{\hfill $\diamond$\endlist}

\usepackage{graphicx} 

\hypersetup{
	colorlinks=true,
	linktocpage=true,
	linkcolor=[RGB]{161,1,49},
	citecolor=[RGB]{10,136,169},      
	urlcolor=cyan,
}
\newtheorem{theorem}{Theorem}
\newtheorem*{main theorem}{Theorem \ref{thm:main-gaussian}}

\newtheorem{corollary}[theorem]{Corollary}
\newtheorem{lemma}[theorem]{Lemma}
\newtheorem{proposition}[theorem]{Proposition}

\newtheorem{definition}[theorem]{Definition}

\theoremstyle{definition}
\newtheorem{remark}[theorem]{Remark}

\title{Cobham's Theorem for the Gaussian integers}
\author{Álvaro Bustos-Gajardo}
\address{Facultad de Matem\'aticas, Pontificia Universidad Cat\'olica de Chile, and Departamento de Matem\'aticas, Facultad de Ciencias, Universidad de Chile}
\email{alvaro.bustos@uchile.cl}
\thanks{The first author was supported by ANID/FONDECYT Postdoctorado 3230159 (year 2023).}
\author{Robbert Fokkink}
\address{Department of Applied Mathematics, Delft University of Technology}
\email{r.j.fokkink@tudelft.nl}
\author{Reem Yassawi}
\address{School of Mathematical Sciences, Queen Mary University of London}
\email{r.yassawi@qmul.ac.uk}
\date{\today}
\keywords{Automata sequences; Gaussian numeration; Cobham's theorem}
\subjclass[2020]{11B85, 13F07, 68Q45}
\thanks{The third author was supported by  the EPSRC, grant number EP/V007459/2.}

\begin{document}

\begin{abstract}
Assuming the four exponentials conjecture, Hansel and Safer showed that if a  subset $S$ of the Gaussian integers  is both $\alpha=-m+\ii $- and $\beta=-n+\ii $-recognizable, then it is syndetic, and they conjectured that $S$ must be eventually periodic. Without assuming the four exponentials conjecture, we show that if $\alpha$ and $\beta$ are multiplicatively independent Gaussian integers, and at least one of $\alpha$,  $\beta$ is not an $n$-th root of an integer, then any $\alpha$- and $\beta$-automatic  configuration  is eventually periodic; in particular we prove Hansel and Safer's conjecture.  Otherwise, there exist non-eventually periodic configurations which are $\alpha$-automatic  for any  root of an integer $\alpha$.    Our work generalises the  Cobham-Semenov theorem to Gaussian numerations. 
\end{abstract}

\maketitle

\section{Introduction}

A sequence $a = (a_n)_{n \geq 0}$ is  {\em $k$-automatic} if $a_n$ is the output
 of a function computed by a finite state automaton, taking as input 
 the base-$k$ expansion of~$n$. 
Such sequences are the simplest sequences beyond the eventually periodic sequences. The base of the numeration system matters:
If $k$ and $\ell$  are multiplicatively independent positive integers, then Cobham's theorem  \cite{Cobham-1969} tells us that  a sequence  is both $k$- and $\ell$-automatic if and only if it is eventually periodic. This theorem is arguably the most significant result 
concerning automatic sequences; its importance is reflected in the fact that it is connected to many deep generalisations and proofs in  different areas of mathematics, such as the work of Semenov and Muchnik in logic \cite{Semenov, Muchnik},  a recent generalisation by Hieronymi and Schulz in \cite{HS-2022 },  the work of Durand and others in dynamical systems and substitution theory \cite{Durand-2008 , Durand-2011, Coven-Dykstra-Lemasurier, Hosseini-Yassawi}, the work of Bès in numeration \cite{Bes-2000}, the work of Adamczewski and Bell in algebra \cite{Bell-2007, Adamczewski-Bell-2008}, and the work of the latter and also  Schäfke and Singer in analysis and ODEs \cite{Adamczewski-Bell-2011, SS-2019}. There are also recent exciting algebraic generalisations of Cobham's theorem by Adamczewski and Faverjon \cite{Adam-Fav-1, Adam-Fav-2}.

While Cobham’s theorem has been extended in many directions, one notable gap is its extension from the integers to Euclidean domains.
In this paper, we take a first step and extend Cobham's theorem to the Gaussian integers.
To formulate the theorem in this context, however, one must first establish what a $\beta$-automatic configuration is for a Gaussian integer $\beta$. In particular, this requires a numeration system with base~$\beta$. 
The numeration system for the Gaussian integers with base  $\beta=-1+\ii $ and digits 0,1 is described by Knuth \cite[Section 4.1]{Knuth:1988}, who attributes it to Penney\cite{Penney}. The only Gaussian integers for which any  $z\in \Z[\ii]$ has an expansion base-$\alpha$   with natural numbers as digits, are  Gaussian integers $\alpha = -m\pm\ii$, with digit set  $\{0,1,\dots , m^2\}$, as shown by Kátai and Szabó in \cite{Katai-Szabo}.

 In \cite{Hansel-Safer}, Hansel and Safer considered numeration systems with nonnegative integer digit sets, and with base given by a Gaussian integer  $\alpha = -m+\ii$  where $m\geq 1$.
  Any two such Gaussian integers  must be multiplicatively independent, see, eg, Theorem~\ref{thm:mult-ind}. Assuming the four exponentials conjecture, they showed that  for such  Gaussian integers,  if the characteristic function  $(a_z)_{z\in \mathbb Z[\ii]}$    of an infinite $S\subset \Z[i]$ is both $\alpha$- and $\beta$-recognizable, then it is syndetic. They  conjectured that this set must be eventually periodic, i.e.,  there exist $\Z$-linearly independent $p,q\in \Z[\ii]$ such that $a_z=a_{z+p}=a_{z+q}$ whenever $z\in \Z[\ii] \backslash F$ for some finite set $F$.

As a result, this extension of Cobham was left unsettled because it was thought to depend on deep analytic number theory. This paper shows that it does not, surprisingly, and opens up the way to Cobham for general number rings. We prove Cobham's theorem for the largest possible family of Gaussian numerations, with any generating finite digit set and without using the four exponentials conjecture. We show
\begin{theorem}\label{thm:main-gaussian}
Let $\alpha$ and $\beta$ be two multiplicatively independent Gaussian integers with $|\alpha|, |\beta|>1$. 
If one of $\alpha$ or
 $\beta$ is not the root of an integer,
then a configuration $(a_z)_{z\in \mathbb Z[\ii]}$ is  $\alpha$- and $\beta$-automatic if and only if  it is eventually periodic. 
\end{theorem}

Niven's theorem limits Gaussian integers that are roots of integers  to have an argument that is a multiple of $\pi/4$ (Theorem~\ref{thm:niven}).
In particular, Hansel and Safer's conjecture is proved true, since $-m+\ii $ is not the root of an integer if $m>1$. 
 The requirement that one of $\alpha$ and $\beta$ must be a  root of an integer cannot be relaxed, as the characteristic function of $\mathbb N$  can be shown to be $\alpha$-automatic whenever $\alpha$ is a root of an integer; see Remark~\ref{rem:best}. Multiplicatively independent Gaussian roots of integers will however obey a Semenov-type theorem \cite{Semenov,Muchnik}, see also the exposition \cite{BHMV1}, behaving like higher dimensional automatic configurations generated by a base-$k$ numeration system for $\mathbb N\times \mathbb N$, where projections of  $\alpha$- and $\beta$-automatic configurations are eventually periodic along any horizontal or vertical  direction.

We base our proof  on Thijmen Krebs' alternative proof of Cobham's theorem\cite{Krebs}. In Section~\ref{sec:preliminaries}, we provide the required background on Gaussian numerations, automaticity and Dirichlet approximation in $\Z[\ii]$. In Section~\ref{sec:periods-1}, we use automaticity and the pumping lemma to generate 2 linearly independent periods, provided that one of  the bases is not a root of an integer. In Section~\ref{sec:periods}, we discuss two closely related notions of periodicity that we will use, and a version of Fine and Wilf's theorem that will guarantee that local periods can propagate to global periods. Finally in Section~\ref{sec:main} we prove our main result.


\section{Preliminaries}\label{sec:preliminaries}

\subsection{Gaussian numerations} In what follows, we discuss some basic facts about the Euclidean domain of Gaussian integers and their numeration systems. For a more in-depth study of Gaussian numeration systems, we direct the reader to \cite[Chapter 7]{Lothaire}.
For $z\in \mathbb Z[\ii]$, we let $|z|$ denote the Euclidean norm of $z$. 
    Given a Gaussian integer $\gamma$ with  $\lvert\gamma\rvert>1$ and a finite subset $D\subset\Z[\ii]$ containing $0$, to each word $w=w_{n-1}\dotsc w_0\in D^*$ we associate the number
        \[[w]_\gamma = \sum_{j=0}^{n-1} w_j\gamma^j; \]
        note that $w_0$ is the least significant digit.
   For $L\subseteq D^*$ we write $[L]_\gamma = \{[w]_\gamma\::\:w\in L \}$. We say that $w\in D^*$ is a $(\gamma,D)$-\emph{expansion} of  $z\in\Z[\ii]$ if $[w]_\gamma=z$.  We say that $(\gamma,D)$ is an {\em integral numeration system} 
 if every Gaussian integer $z$ has a unique  expansion  $n=  \sum_{j=0}^{k} d_j\gamma^j$ with $d_j\in D$ and $d_{k }\neq 0$. In this case we write 
  $(n)_\gamma = d_k \dots d_0$. A set $D\in \Z[\ii]$ is a {\em complete residue system} for $\Z[\ii] \pmod \gamma$
 if every Gaussian $z$ is congruent modulo $\gamma$ to a unique
element of $D$. 
We record the following result, which indicates that we have an abundance of Gaussian numeration systems.
\begin{theorem}\label{thm:numerations-exist}
Let $0\neq \gamma=a+bi\in \Z[\ii]$, let $N\coloneq |\gamma|^2$
and let $\lambda = \gcd(a,b)$.  Suppose  that $|\gamma|>1$. Then:    \begin{enumerate}\item\cite[Theorem 7.4.1]{Lothaire} The set
\[ D(\gamma)=\biggl\{p+\ii q: p= 0,1, \dots , \frac{N}{\lambda} -1, q=0,1, \dots , \lambda -1 \biggr\}\]
is a complete residue system for $\Z[\ii] \bmod \gamma$.
\item\cite[Proposition 7.4.2]{Lothaire} Let $D$ be a finite set in $\Z[\ii]$. If every Gaussian integer has an expansion  using $(\gamma ,D)$, then $(\gamma, D)$ is a numeration system if and only if $D$ is a complete residue system for $\Z[\ii] \mod \gamma$ and $0\in D$.
\item \cite{DDG}  Let $|\gamma|\geq \sqrt{5} $ and let
\[D= \biggl\{d \in \Z[\ii] :  \operatorname{Re}\Bigl(\frac{d}{\gamma}\Bigr),\operatorname{Im}\Bigl(\frac{d}\gamma{}\Bigr)\in \Bigl[-\frac{1}{2}, \frac{1}{2}\Bigr)\biggr\}.\]

Then $(\gamma, D)$ is a numeration system.
\item\cite{Katai-Szabo} Let 
$D=\{0,1, \dots ,N-1\}.$
Then 
$(\gamma, D) $ is a numeration system for $\Z[\ii]$ if and only if $\gamma = -n\pm \ii$. 
\end{enumerate}
\end{theorem}

\subsection{\texorpdfstring{Dirichlet approximation in $\Z[\ii]$}{Dirichlet approximation in Z[\ii ]}}
Recall that $\alpha, \beta\in \mathbb{Z}[\ii]$ are {\em multiplicatively independent} if there are no natural numbers $j,k$ such that $\beta^j= \alpha^k$.
In his original proof over $\mathbb N$, Cobham needed the fact that for multiplicatively independent numbers  $\alpha, \beta\geq 2$,
 the multiplicative group generated by $\alpha$ and $\beta$ is dense in the positive real line. It turns out that this result is very hard to generalise to $\mathbb C$.  Hansel and Safer showed that if $\alpha$ and $\beta$ are multiplicatively independent Gaussians, then the four exponentials conjecture implies that the multiplicative group
 $\langle\alpha,\beta\rangle^\times=\{\alpha^m\beta^n\::\:m,n\in\Z\}$
   is dense in  $\mathbb C$. In our proof of Theorem~\ref{thm:main-gaussian}, we only need that there are no isolated points in  $\langle\alpha,\beta\rangle^\times$, a much weaker result.  

\begin{lemma}[\hspace{1sp}\cite{Bosma-Fokkink-Krebs}, Lemma~3] Let $\alpha,\beta\in\Z[\ii]$ be two multiplicatively independent Gaussian integers with $\lvert\alpha\rvert,\lvert\beta\rvert>1$. Then $1$ is not an isolated point in the multiplicative group $\langle\alpha,\beta\rangle^\times$.
\end{lemma}
\begin{proof}
     We need to prove that $\alpha$ and
    $\beta$ are multiplicatively dependent if $1$ is isolated in $G=\langle\alpha,\beta\rangle^\times$. 
    Suppose that $d\in \mathbb C\setminus\{0\}$ is a density point of $G$, say that it is the limit of $g_n\in G$. 
    Then $g_{n+1}/g_{n}$ converges to $1$, and since $1$ is isolated this means that the sequence
    is eventually constant. It follows that $G$ is a discrete subset of the punctured plane if $1$ is isolated.
    Let $H\subset G$ be the cyclic subgroup that is generated by $\alpha$. Every residue class of the quotient $G/H$ 
    has a representative in the annulus $A=\{z : 1 \leq |z| \leq |\alpha|\}$. Since $G\cap A$ is compact and discrete,
    it is finite. Therefore, $G/H$ is finite and $\beta^n\in H$ for some $n$, which means that $\alpha$ and $\beta$
    are multiplicatively dependent.
\end{proof}

\begin{corollary}\label{cor:dirichlet}
    Let $\varepsilon>0$. Given two non-zero multiplicatively independent Gaussian integers $\alpha,\beta\in\Z[\ii]$ with $\lvert\alpha\rvert,\lvert\beta\rvert>1$, there exist $m,n>0$ such that
    \begin{equation}\label{eq:relatively-close-powers}
        \lvert \alpha^m-\beta^n\rvert<\varepsilon\lvert \beta^n\rvert.
    \end{equation}
\end{corollary}

\begin{proof}
    Without loss of generality, we may assume that $\lvert\alpha\rvert>\lvert\beta\rvert$, as the only Gaussian integers of absolute value $1$ are the units, and thus two Gaussian integers of the same absolute value must be multiplicatively dependent. Since $1$ is not an isolated point in $\langle\alpha,\beta\rangle^\times$, there exist infinitely many values $m,n\in\Z$ for which $\lvert \alpha ^m \beta^{-n} - 1\rvert<\varepsilon$. Furthermore, $m$ and $n$ must have the same sign for all sufficiently small $\varepsilon$, since otherwise $\lvert \alpha^m \beta^{-n} - 1\rvert$ is bounded from below by $\min(\lvert \alpha\beta\rvert - 1, 1- \lvert \alpha\beta\rvert^{-1})>0$.
    
    Once again, without loss of generality, we may assume for sufficiently small $\varepsilon$ that $m>0$ due to the continuity of the function $z^{-1}$ near $1$. Thus, by multiplying our inequality by $\lvert \beta^n\rvert$, we obtain \eqref{eq:relatively-close-powers}.
\end{proof}

The following classical result
appears as Corollary 3.12 in \cite{niven} and is known as
Niven's theorem.
\begin{theorem}\label{thm:niven}
    Let $\alpha\in\mathbb{Z}[\ii ]$ be such that $\alpha^j\in \mathbb Z$
    for some $j\in\mathbb N$. Then $\text{arg}(z)$ is a multiple of~$\pi/4$.
\end{theorem}
 It takes some effort to decide whether two numbers in $\mathbb{Z}[\ii ]$   are multiplicatively dependent.  Hansel and Safer considered bases $\alpha=-m+\ii $ for $m\geq 1$. It turns out that
the only Gaussian integers that are multiplicatively dependent on $\alpha$ are powers
of $\alpha$.
\begin{theorem}\label{thm:mult-ind}
    Let $\alpha=-m+\ii $ for $m\geq 1$. If $\beta^j=\alpha^k$ for $j,k\in\mathbb N$,
    then, up to a root of unity, $\beta$ is a power of $\alpha$.
\end{theorem}
\begin{proof}
    Let $y$ denote the norm of $\beta$; then $y^j=(m^2+1)^k$.
    Therefore, $j\cdot\text{ord}_p(y)=k\cdot\text{ord}_p(m^2+1)$ for all (rational) primes $p$. 
    We first treat the case that $\text{gcd}(j,k)=1$.
    Then $k$ divides $\text{ord}_p(y)$. In particular, $y=z^k$ for
    some natural number $z$. It follows that $z^{jk}=(m^2+1)^k$ and we conclude that $z^j=m^2+1$.
    By
    Mih{\u{a}}ilescu's theorem~\cite{mihailescu2004}, $j=1$ and $\beta=\alpha^k$. In the general case $j$
    and $k$ have a common factor $h$, it follows that $\beta^h=\alpha^{hk}$. Up to a root of unity,
    $\beta$ is equal to $\alpha^k$.
\end{proof}
Even the seemingly simple task of deciding whether 
$\alpha$ and $\beta$ are multiplicatively dependent runs into a weak 
form of Catalan's conjecture, wherein one of the powers is a square.
Such number theoretic difficulties naturally arise when considering
the extension of Cobham's theorem from the integers to other domains.
This may explain why this extension has been neglected so far.

\subsection{Finite automata and automaticity}

    A \emph{deterministic finite automaton with output} (DFAO) is defined by a $6$-tuple $\mathcal{A}=(S,D,\delta,s_0,A,\tau)$, where $S$ is a finite set of \emph{states}, $s_0\in S$ is the \emph{initial state}, $D$ and $A$ are finite sets called the \emph{input} and \emph{output alphabet}, respectively, $\delta\colon S\times D\to S$ is the \emph{transition function} and $\tau\colon S\to A$ is the \emph{output function}. The transition function extends to the set $D^*$ of all words with letters in $S$ via the recurrence relation:
        \[\delta(s,wa)=\delta(\delta(s,w),a),\quad w\in D^*,\,a\in D, \]
    and thus every DFAO determines a function $D^*\to A$ given by $w\mapsto\tau(\delta(s_0,w))$.
We shall use DFAOs by feeding in corresponding expansions of Gaussian integers in terms of powers of a fixed $\gamma$ such that $D$ contains a complete residue system for $\Z[\ii] \bmod \gamma$; when we want to emphasize this connection between the automaton and the associated $\gamma\in\Z[\ii]$, we include appropriate subscripts in the corresponding sets, say, $\mathcal{A}_\gamma=(S_\gamma,D_\gamma,\delta_\gamma,s_{0,\gamma},A_\gamma,\tau_\gamma)$.

We will need to consider automatic configurations defined by DFAOs with redundant digit sets, i.e. Gaussian numerations where each Gaussian integer has at least one expansion with nonzero leading digit.

\begin{definition}\label{def:automaton}
    Let   $\gamma\in\Z[\ii]$ be a  Gaussian integer with $|\gamma|>1$ and let $D\subset\Z[\ii]$ be a finite subset   for which every $z\in\Z[\ii]$ has  at least one finite expansion with digits from $D$. Let  $\mathcal{A}=(S,D,\delta,s_0,A,\tau)$ be a DFAO. A configuration $(a_z)_{z\in \Z[\ii]}$ is $(\gamma,D)$-\emph{automatic} if  
   $a_z=\tau(\delta(s_0, w))$ whenever $[w]_\gamma =z$. 
\end{definition}

\begin{figure}
	\center{\includegraphics[scale=.25]{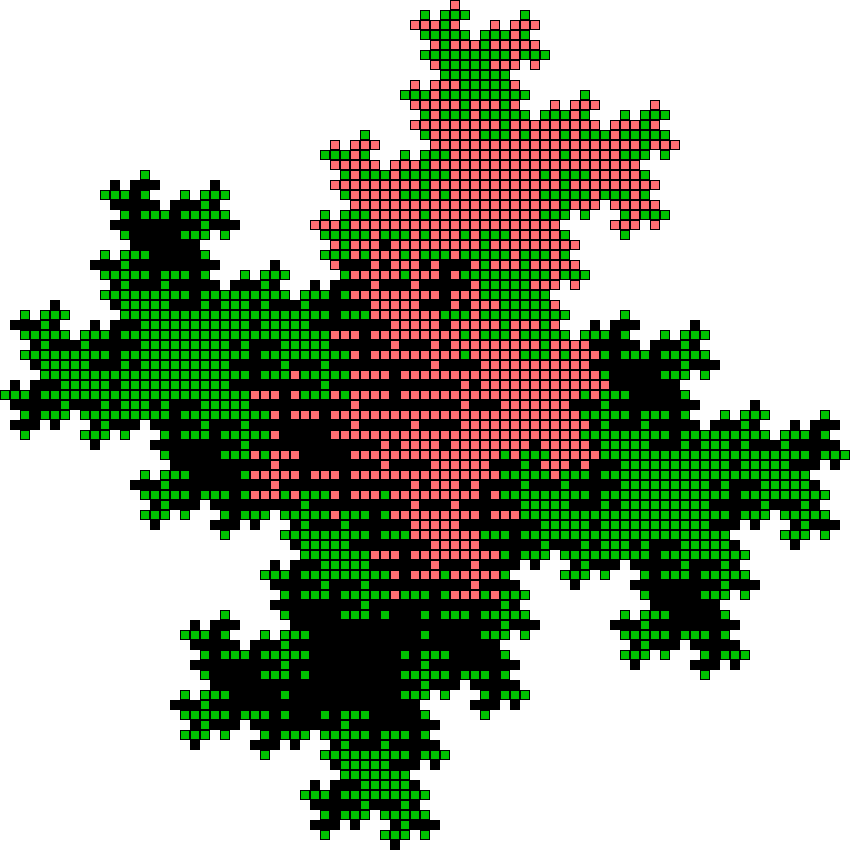} \begin{tikzpicture}[scale=0.8]
	\node[draw=red,line width=0.5mm,shape=circle] (qA) at (150:2cm) {\color{red}$A$};
	\node[draw=black,line width=0.5mm,shape=circle] (qB) at (270:2cm) {$B$};
	\node[draw=green!50!black,line width=0.5mm,shape=circle] (qC) at (30:2cm) {\color{green!50!black}$C$};
	\draw[-latex,looseness=6] (qA) edge node [above] {\small$0,1$} (qA);
	\draw[-latex,looseness=6] (qC) edge node [above] {\small$1,-1$} (qC);
	\draw[-latex,bend left=15] (qA) edge node [above] {\small$\ii,-\ii$} (qC);
	\draw[-latex,bend left=15] (qC) edge node [below] {\small$\ii,-\ii$}  (qA);
	\draw[-latex,looseness=8,rotate=180] (qB) edge node [below] {\small$0,-1$} (qB);
	\draw[-latex,bend right=15] (qA) edge node [left] {\small$-1$} (qB);
	\draw[-latex,bend right=15] (qB) edge node [right] {\small$1$} (qA);
	\draw[-latex,bend right=15] (qC) edge node [left] {\small$0$} (qB);
	\draw[-latex,bend right=15] (qB) edge node [right] {\small{$\ii,-\ii$}} (qC);
	\draw[latex-] (qA) -- +(-1,0);
    \draw[color=white] (0,-5);
\end{tikzpicture}}
	\caption{A $(1+2\ii)$-automatic configuration with digit set $D=\{0,\pm 1,\pm \ii\}$. The corresponding automaton is shown to the right; the colours correspond to its states, which we also take to be the output alphabet.} 
	\label{fig:automatic-1}
\end{figure}

See Figure~\ref{fig:automatic-1} for an example of a $(1+2\ii, D)$-automatic configuration with  $D=\{0,\pm 1,\pm \ii\}$.
Examples of $\beta$-automatic configurations include characteristic functions of  $\beta$-recognizable sets, also called $\beta$-automatic sets. Any automatic configuration  can be seen as a  partition of $\Z[\ii]$ into automatic sets.

In this paper, almost all DFAOs are direct reading, that is, a word which represents $z$ is fed into the machine starting with its most significant digit. One exception is when we use Lemma~\ref{lem:transducer}. This will not be an issue as a configuration is $(\gamma,D)$-automatic in direct reading if and only if it is $(\gamma,D)$-automatic in reverse reading.

  In what follows, we shall assume that $\alpha,\beta\in\Z[\ii]$ are two multiplicatively independent Gaussian integers, and we shall use the letter $\gamma$ to refer generically to either of them.

\begin{remark}
     We  say that a DFAO $\mathcal{A}_\gamma$ is \emph{consistent} if, whenever $[w_1]_\gamma=[w_2]_\gamma$, then $\tau_\gamma(\delta_\gamma(s_{0,\gamma},w_1))=\tau_\gamma(\delta_\gamma(s_{0,\gamma},w_2))$. If $(\gamma, D)$ is integral, then any $(\gamma,D)$-automaton is consistent. Note that  Definition~\ref{def:automaton} implicitly requires that the DFAO defining  a $(\gamma , D)$-automatic configuration is consistent. Henceforth all DFAOs we work with are assumed consistent. 
\end{remark}

The classical definition of automaticity  uses an integral numeration system $(\gamma, D)$.  Proposition~\ref{prop:automaticity-equivalent} below tells us that if 
$(\gamma, D)$ is integral and $D\subset D'$,
then the family of $(\gamma, D)$-automatic configurations   equals  the family of $(\gamma, D')$-automatic configurations generated by consistent automata. For this reason, if $(a_z)_{z\in \Z[\ii]}$ is $(\gamma ,D)$-automatic, then we will sometimes write that $(a_z)_{z\in \Z[\ii]}$ is $\gamma$-automatic if we do not need to state what the digits are.

 The following is based on \cite[Lemma 7.1.1]{Lothaire};
we include it for completeness. It tells us how to transform a non-integral expansion  to an integral one. The required formal definitions are in the body of the proof.
\begin{lemma}\label{lem:transducer}Let $|\beta|>1$,
    let $(\beta, D)$ be an integral Gaussian numeration and let $D\subset D'$. Then there exists a 1-uniform transducer with a terminal function that sends 
    any $(\beta, D')$ expansion of $z$  to $(z)_{\beta,D}\in D^*$.
\end{lemma}
\begin{proof}
The idea is as follows. We go through each digit of an expansion $d_n'\dotsc d_0'\in (D')^*$ for a Gaussian integer $z\in\Z[\ii]$, starting from the least significant digit $d_0'$, which we replace by the unique digit $d_0\in D$ which is congruent to $d_0'$ modulo $\beta$. We now have a carry $c_0$ for which the equality $\beta c_0=d_0'-d_0$ holds, so $z=[d_n'\cdots d_1'd_0]_\beta + \beta c_0$. We iterate this process, at each step defining $c_k=(d_k'-d_k)/\beta$ and $d_{k+1}$ as the only element of $D$ congruent to $c_{k}+d_{k+1}'$  modulo $\beta$, so that $d_{k+1}'+c_k=\beta c_{k+1} + d_{k+1}$ for some value $c_{k+1}\in\Z[\ii]$, which will become the next carry. At each step the equality $z=[d_n'\dotsc d_{k+1}'d_k\dotsc d_0]_\beta + \beta^{k+1} c_k$ holds. After we modify all of the digits $d_0',\dotsc,d_n'$, we have an expression of the form $z= \beta^{n+1} c_n+[d_n\dotsc d_0]_\beta$, so we only need to append the digits of the $(\beta,D)$-expansion of $c_n$ at the start to get a full $(\beta,D)$-expansion of $z$.

For this procedure to work properly, at each step we need to keep track of the value of the carry $c_k$. It can be seen that, if we start with an initial carry of $0$, then there are finitely many choices for each $c_k$, and so this process can be carried out by a transducer whose set of states $Q$ corresponds to all possible carries. Indeed, set $m=\max\{\lvert d'-d\rvert : d'\in D',d\in D\}$ and $c=m/(|\beta|-1)$. Note that if $|c_k|<c$, then $|c_{k+1}|\leq (|c_k|+|d-d'|)/|\beta|< \left|(c+m)/\beta \right|= c$. This shows that  $Q$ is finite.


Formally, consider the subsequential finite transducer $\mathcal A=(Q,D',\delta, q_{0}, D, \tau)$
  defined as follows: $Q= \{s\in \Z[\ii]:|s|<c\}$ is the set of possible carries, and the set of edges, which defines 
  $\delta:Q\times D'\rightarrow Q$ and $\tau:D'\rightarrow D$, is
\[ \{s\stackrel{d'/d}{\longrightarrow}s': s+d'=\beta s'+d\}.\]
Let $d_n'\dots d_0'\in D'^*$ be a $( \beta, D')$-expansion of $z$. Starting at $s_0=0$ there is a unique path
\[s_0\stackrel{d'_0/d_0}{\longrightarrow}s_1 
\stackrel{d'_1/d_1}{\longrightarrow}s_2
\stackrel{d'_2/d_2}{\longrightarrow}
\dots
\stackrel{d'_{n-1}/d_{n-1}}{\longrightarrow}s_n
\stackrel{d'_n/d_n}{\longrightarrow}s_{n+1}.
\]
Define the terminal function $t:Q\rightarrow D^*$ as $t(s)=(s)_{(\beta,D)}$.
It can be verified that 
$z= \sum_{j=0}^n d_j\beta^j
+s_{n+1}\beta^{n+1}$, hence the output of the transducer on inputting $d_n'\dots d_0'\in D'^*$ is \[t(s_{n+1}) d_n\dots d_0  = (z)_{(\beta, D)}.\]\end{proof}

The following proposition is reminiscent of Cabezas and Leroy's \cite[Proposition 4.1]{Cabezas-Leroy}, which gives a similar result in the primitive, aperiodic, substitutional case.
\begin{proposition}\label{prop:automaticity-equivalent}
Let $|\beta|>1$,
    let $(\beta , D)$ be an integral Gaussian numeration system and let $D\subset D'$. Then $(a_z)_{z\in \Z[\ii]}$ is $(\beta , D)$-automatic if and only if it is 
    $(\beta, D')$-automatic and generated by a consistent automaton. 
\end{proposition}
\begin{proof}
Let $(Q,D',\delta,q_{0,T},D, \tau, t)$ be the transducer with terminal function $t$ as defined in Lemma~\ref{lem:transducer}, and let $(Q_D, D, \delta_D,q_{0,D}, A_D,\tau_D)$ be the $(\beta , D)$-automaton that generates the $(\beta, D)$-automatic configuration $(a_z)_{z\in \Z[\ii]}$. 
We construct a $(\beta, D')$-automaton
$\mathcal A=(Q\times Q_D,D',\delta_{D'}, [q_{0,D},q_{0,T}], A_D, \tau_{D'}   )$ 
by defining  $\delta_{D'}:Q\times Q_D \times D'\rightarrow Q\times Q_D$ and $\tau_{D'} :Q_D\rightarrow A_D$ as
\[ \delta_{D'}([p,q],d') = [\delta(p,d'), \delta_D(q,\tau(d'))] \mbox{ and }   \tau_{D'} ([p,q]) = \tau_D(\delta_D(q,t(p))).\] 
One uses induction to verify that $\mathcal A$ generates $(a_z)$ in  $(\beta, D')$-reading. 
Note furthermore that the automaton  $\mathcal A$ constructed in Proposition~\ref{prop:automaticity-equivalent} is consistent.\\
Conversely, if $\mathcal A$ is a consistent automaton generating the $(\beta, D')$-automatic configuration $(a_z)_{z\in\Z[\ii ]}$, then one can restrict to $D$-inputs to see that the configuration is also $(\beta, D)$-automatic.
\end{proof}

Let $B(z,r)\subset \C$ denote the  ball centred at $z$ with radius $r$. Next we  choose an appropriate set $D_\gamma$, which by the next lemma will be  large enough to ensure that every Gaussian integer relatively close to $0$ has a short $(\gamma, D_\gamma)$-expansion.

\begin{lemma}\label{lem:small-numbers-have-short-expansions}
     Let $(\gamma, D)$ be an integral numeration  system. Let $r\geq 2$ be chosen so that 
     \begin{equation} \label{eq:expanded_digit_set}
        D\subset D_\gamma \coloneq B(0,r\lvert\gamma\rvert)\cap\Z[\ii].
    \end{equation}
If  $z\in\Z[\ii]$ satisfies
     $\lvert z\rvert \le r\lvert\gamma\rvert^n$, then  $z$ has at least one $(\gamma, D_\gamma)$-expansion of length at most $n$.
\end{lemma}
\begin{proof}
We need to show that if  $\lvert z\rvert \le r\lvert\gamma\rvert^n$, then there exists a word $w=d_{n-1}\dotsc d_{0}\in D_\gamma^n$ for which $[w]_\gamma = z$.
    We proceed by induction on $n$; the result is evidently true if $\lvert z\rvert\le r\lvert\gamma\rvert$ by our choice of $D_\gamma$. Suppose that the result is valid for all Gaussian integers in $ B(0,r\lvert \gamma\rvert^n)\cap\Z[\ii]$, and let $z\in B(0,r\lvert\gamma\rvert^{n+1})\cap\Z[\ii]$. We have $\lvert\frac{z}{\gamma}\rvert\le \frac{r\lvert\gamma\rvert^{n+1}}{\lvert\gamma\rvert}=r\lvert\gamma\rvert^n$; we claim that this implies that there exists a Gaussian integer $q\in\Z[\ii]$ with $|q|\leq r\lvert\gamma\rvert^n$ and
    $|q-\frac{z}{\gamma}|\leq 
    \sqrt{2}$.  Suppose first that $\frac{z}{\gamma}=u_1+u_2\ii$ with $u_1,u_2\in\R_{\ge 0}$ and set $q=\lfloor u_1\rfloor + \lfloor u_2\rfloor\ii$. Then, $\lvert q\rvert^2=\lfloor u_1\rfloor^2+\lfloor u_2\rfloor^2\le u_1^2+u_2^2\le (r\lvert\gamma\rvert^n)^2$.  The proof for the other three quadrants proceeds  
    in the same fashion, replacing $\lfloor u_j\rfloor$ by $-\lfloor -u_j\rfloor$ when that coordinate is negative. The inequality $\lvert q - \frac{z}{\gamma}\rvert \le\sqrt{2}$  follows immediately from the fact that $\lvert x-\lfloor x\rfloor\rvert <1$.

    As $\lvert q\rvert\le r\lvert\gamma\rvert^n$, by the induction hypothesis there must exist $d_0,\dotsc,d_{n-1}\in D_\gamma$ such that $q=d_0+\cdots+d_{n-1}\gamma^{n-1}$, i.e., the word $d_{n-1}\dotsc d_0$ is a $(\gamma , D_\gamma)$-expansion of $q$ of length at most $n$. Thus, $q\gamma=d_0\gamma+\cdots+d_{n-1}\gamma^n$. Furthermore, since $\lvert z- q\gamma\rvert = \lvert \gamma \rvert \cdot \lvert \frac{z}{\gamma}-q\rvert  \le \lvert\gamma\rvert\cdot\sqrt{2}< r\lvert\gamma\rvert$, the Gaussian integer $d^*=z-q\gamma$ is in the ball $B(0,r\lvert\gamma\rvert)$ and thus belongs to $D_\gamma$. Hence
        \[z= q\gamma + (z-q\gamma)=d_{n-1}\gamma^n+\cdots+d_0\gamma+d^*, \]
    where all digits $d_0,\dotsc,d_{n-1},d^*\in D_\gamma$, so that  the word $d_{n-1}\dotsc d_0d^*\in D_\gamma^{n+1}$ is a $(\gamma, D_\gamma)$-expansion of $z$ of length
    at most $n+1$.
\end{proof}
%
%
%
%

\section{Generating periods}\label{sec:periods-1}

In this section, in Lemma~\ref{lem:shifting-expansions} we exploit the structure of $(\gamma, D)$-expansions of Gaussian integers to show repetitive behaviour on a corresponding automatic configuration. Lemmas~\ref{lem:non-collinearity-1}~and~\ref{lem:non-collinearity-with-two-automata} will allow us to later create pairs of non-collinear vectors which we will bootstrap, to show that configurations satisfying the conditions of Theorem~\ref{thm:main-gaussian} are eventually periodic.

Let $\gamma$ be a Gaussian integer with $\lvert\gamma\rvert>1$, and let $D$ be a digit set for which each $z\in\Z[\ii]$ has at least one $(\gamma , D)$-expansion. Suppose that $\mathcal{A}_\gamma=(S,D,\delta,s_{0},A,\tau)$ is a consistent $(\gamma , D)$-automaton, and let $s\in S$. We define
    \begin{equation}\label{eq:words-ending-at-s}
        L_{\gamma,s}\coloneqq \{w\in D^*\::\: \delta(s_{0},w)=s\}.
    \end{equation}
The set $[L_{\gamma,s}]_\gamma$ corresponds to all the Gaussian integers which have at least one $(\gamma , D)$-expansion that ends at state $s$ in the automaton $\mathcal{A}_\gamma$. The sets $\{[L_{\gamma,s}]_\gamma\::\:s \in S\}$ are a cover of $\Z[\ii]$, but they might not be a partition as some Gaussian integers might have more than one $(\gamma, D)$-expansion.

\begin{lemma}\label{lem:shifting-expansions}
    Let $\mathcal{A}=(S,D,\delta,s_{0},A,\tau)$ be a $(\gamma , D)$-consistent DFAO, and let $(a_z)_{z\in\Z[\ii]}$ be the corresponding automatic configuration. Then, for every $x,y\in [L_{\gamma,s}]_\gamma$, any $n\in\N$ and every $z\in[D^n]_\gamma$, we have
        \begin{equation}\label{eq:shifting-expansions-formula}
            a_{x\gamma^n+z}=a_{y\gamma^n+z}.
        \end{equation}
   \end{lemma}

\begin{proof}
    By definition, there exist $w_1,w_2\in L_{\gamma,s}$ such that $[w_1]_\gamma=x,[w_2]_\gamma=y$. Thus, for any $v\in D^n$, we have $\delta(s_0,w_jv)=\delta(\delta(s_0,w_j),v)=\delta(s,v)$. Hence, as $[w_jv]_\gamma=[w_j]_\gamma\gamma^n+[v]_\gamma,j\in\{1,2\}$, if we write $z=[v]_\gamma$, we have
    \[a_{x\gamma^n+z}=\tau(\delta(s_0,w_1v))=\tau(\delta(s,v))=\tau(\delta(s_0,w_2v))=a_{y\gamma^n+z},\]
    since $w_1v$ and $w_2v$ are $(\gamma, D)$-expansions for $x\gamma^n+z$ and $y\gamma^n+z$, respectively.
\end{proof}

In particular, for $D=D_\gamma=\Z[\ii]\cap B(0,r\lvert\gamma\rvert)$ as defined in \eqref{eq:expanded_digit_set}, the above lemma holds for any $z$ with $\lvert z\rvert\le r\lvert \gamma^n\rvert$, as per Lemma~\ref{lem:small-numbers-have-short-expansions}. Henceforth when we use the notation $D_\gamma$ we mean an enlarged digit set satisfying
\eqref{eq:expanded_digit_set}.

 Let $S_\infty$ be the set of states in $S_\beta$ for which $[L_{\beta, s}]_\beta$ is infinite, where $L_{\beta, s}$ is defined as in \eqref{eq:words-ending-at-s}. We will need $S_\infty$ to contain certain  triples of non-collinear points. The following lemma tells us that this is guaranteed if $\beta$ is not 
 an $n$-th root of an integer.  
 We use the existence of a cycle rooted at $s \in S_\infty$ to deduce that Gaussian integers with expansions ending at $s$ do not live in a one-dimensional subspace unless $\beta$ is restricted. 
 The existence of such a cycle is in general not guaranteed, but the pumping lemma can be used to create a weaker version of this condition which will be sufficient.
 For a word $w$ over an alphabet $A$ we let $|w|$ denote the length of $w$.
  (It will be clear from the context that  the object considered is a word $w$ and not a Gaussian integer $z$.)

   \begin{definition}\label{def:return-number} Let $\mathcal A_\beta $ be a DFAO, and  let $s\in S_\infty$. Then, by the pigeonhole principle, there exist $s'\in S_\infty$ and three words $w_1,w_2,w_3\in D_\beta^*$ for which
    $\delta_\beta(s_0,w_1)=s'$,  $\delta_\beta(s',w_2)=s'$ and $\delta_\beta(s',w_3)=s$. 
    Furthermore, since $s'\in S_\infty$, there are infinitely many choices for  $[w_1]_\beta$, so we pick one satisfying   $[w_1]_\beta \cdot (1-\beta^{\lvert w_2\rvert}) \ne  [w_2]_\beta$. Note that for each $j$, $u_j=[w_1(w_2)^jw_3]_\beta$ is a number
      whose $(\beta, D_\beta)$-expansion ends at state $s$. 
    We  call such numbers 
      $u_j$ {\em return numbers to $s$}.\end{definition}

     If there is a cycle in the underlying graph of $\mathcal A_\beta$ rooted at $s$, we can take $s=s'$ and $w_3$ to be the empty word in Definition~\ref{def:return-number}.

\begin{lemma}\label{lem:non-collinearity-1}
    Let $s\in S_\infty$, and  suppose that $(u_j)$  are  return numbers to $s$.
      If there exists an  arithmetic progression of indices such that  $u_k,u_{k+\ell},u_{k+2\ell}$ are collinear return numbers, then $\beta$ is the root of an integer.
\end{lemma}

\begin{proof}

We have $u_j=[w_1w_2^j w_3]_\beta$ where the words $w_i$ are as in Definition~\ref{def:return-number}. First suppose that $s=s'$ and $w_3$ is the empty word, i.e., that there is a cycle rooted at $s$.

    Let $N=\lvert w_2\rvert$ and $h=[w_2]_\beta$. The sequence of return numbers $(u_j)_{j\in\N}$  satisfies the recurrence relation $u_{j+1} = \beta^N u_j + h$, from which the following equalities immediately follow:
        \begin{align*}
            u_{k+\ell} &= \beta^{N\ell}u_k+h(\beta^{N(\ell - 1)} + \beta^{N(\ell-2)} + \cdots + \beta ^N + 1)  \\
            u_{k+2\ell} &= \beta^{N\ell}u_{k+\ell}+h(\beta^{N(\ell - 1)} + \beta^{N(\ell-2)} + \cdots + \beta ^N + 1)  \\
            {}\implies u_{k+2\ell}-u_{k+\ell} &= \beta^{N\ell}(u_{k+\ell}-u_k).
        \end{align*}
        
    We now face two scenarios. If $u_{k+\ell}-u_k\ne 0$, we may divide by this term and obtain the equality $\beta^{N\ell}=\frac{u_{k+2\ell}-u_{k+\ell}}{u_{k+\ell}-u_k}$. If $u_k,u_{k+\ell},u_{k+2\ell}$ lie on the same line in the plane, then $u_{k+2\ell}-u_{k+\ell}$ and $u_{k+\ell}-u_{k}$ must lie on a parallel line which passes through the origin, i.e. there exists a real constant $\lambda$ for which $u_{k+2\ell}-u_{k+\ell} = \lambda(u_{k+\ell}-u_{k})$. But then $\beta^{N\ell}=\lambda\in\R$, which is only possible if $\arg(\beta)$ is a rational multiple of $2\pi\ii$.

    Otherwise, if $u_{k+\ell}-u_k=0$, we obtain that $u_{j+\ell}=u_j$ for all $j\ge k$ from the aforementioned recurrence relation, so the sequence $(u_j)_{j\in\N}$ remains bounded. Write $\hat{u}_j=u_j-\frac{h}{1-\beta^N}$; from this definition, it holds that:
        \begin{align*}
            \hat{u}_{j+1}&=u_{j+1}-\frac{h}{1-\beta^N}\\
            &=\beta^N u_j+h-\frac{h}{1-\beta^N} \\
            &=\beta^N\left(\hat{u}_j+\frac{h}{1-\beta^N}\right)+h\left( 1-\frac{1}{1-\beta^N}\right)\\
            &=\beta^N\hat{u}_j+\frac{\beta^Nh}{1-\beta^N} + \frac{h((1-\beta^N)-1)}{1-\beta^N}\\
            &=\beta^N\hat{u}_j,
        \end{align*}
    and since $\lvert\beta\rvert>1$, we have that $\lvert \hat{u}_j\rvert\to\infty$ unless $\hat{u}_0=u_0-\frac{h}{1-\beta^N}=0$. Since $u_0=[w_1]_\beta$ and $h=[w_2]_\beta$, this is impossible as we chose $w_1$ and $w_2$ so that $[w_1]_\beta \cdot (1-\beta^{\lvert w_2\rvert}) \ne  [w_2]_\beta$. As the sequence $(u_j)_{j\in\N}$ is bounded if and only if $(\hat{u}_j)_{j\in\N}$ remains bounded as well, we get a contradiction, so this scenario never happens.
    
    Consider now the general case when $s\neq s'$.   We define $M=\lvert w_3\rvert$, $h' = [w_3]_\beta$, and for each $j$, $u_j=[w_1(w_2)^j]_\beta$, which are return numbers to $s'$.
    Then we set    $v_k=\beta^M u_k+h'$, which are return numbers to $s$, and  each term $v_k$ is obtained from the term $u_k$ via a dilation followed by a translation. As these transformations preserve collinearity, the  return numbers $v_k, v_{k+\ell}, v_{k+2\ell}$ lie on the same line if and only if the respective  $u_k, u_{k+\ell}, u_{k+2\ell}$ are collinear as well, and thus an application of the first case above to the state $s'$ shows us that $\arg(\beta)$ is a rational multiple of $2\pi$. 
\end{proof}


\begin{lemma}\label{lem:non-collinearity-with-two-automata}
   Given an $\alpha$- and $\beta$-automatic configuration, let $s\in S_\infty\subset S_\beta$. If  $\beta$ is not a root of an integer,
     then there exists a state $t=t(s)\in S_\alpha$ for which we may choose three non-collinear Gaussian integers $x_{s,t},y_{s,t},z_{s,t}\in [L_{\beta,s}]_\beta\cap [L_{\alpha,t}]_\alpha$. \end{lemma}
     
     That is, each of the $x_{s,t},y_{s,t},z_{s,t}$ have both a $(\beta, D_\beta)$-expansion ending at state $s$ in the $(\beta , D_\beta)$ automaton, and an $(\alpha, D_\alpha)$-expansion ending at state $t$ in the $(\alpha, D_\alpha)$ automaton.

\begin{proof}
    
    We pick a sequence of return numbers  $(u_j)_{j\in\N}$ to $s$ satisfying the  conditions of Definition~\ref{def:return-number}. 
    Define a partition $\{U_t\}_{t\in S_\alpha}$ of $\N$ satisfying the condition that, whenever $j\in U_t$, the Gaussian integer $u_j$ has an $(\alpha, D_\alpha)$-expansion ending at state $t$ in $\mathcal{A}_\alpha$; such a partition always exists, as every Gaussian integer has at least one $(\alpha, D_\alpha)$-expansion. By van der Waerden's theorem \cite{vdW}, there exists at least one state $t(s)\in S_\alpha$ for which the set $U_{t(s)}$ contains an arithmetic progression $\{k,k+\ell,k+2\ell\}$ of length $3$. By Lemma~\ref{lem:non-collinearity-1}, the trio $u_k,u_{k+\ell},u_{k+2\ell}$ cannot lie on the same line, as otherwise $\arg(\beta)$ would be rational. Furthermore, as all three points lie in $U_{t(s)}$, they each have at least one $(\alpha , D_\alpha)$-expansion ending at state $t(s)$. If we define $x_{s,t(s)}=u_k$, $y_{s,t(s)}=u_{k+\ell}$ and $z_{s,t(s)}=u_{k+2\ell}$, we are done.

\end{proof}

\section{\texorpdfstring{Periodic patterns on $\Z[\ii]$}{Periodic patterns on $Z[\ii ]$}}\label{sec:periods}

We will use the Gaussian integers  $x_{s,t},y_{s,t},z_{s,t}$ from the statement of Lemma~\ref{lem:non-collinearity-with-two-automata} to define vectors which define two-linearly independent periods for an $\alpha$- and $\beta$-automatic configuration. Given two $\Z$-linearly independent
Gaussian integers $p_1$, $p_2$, let $\langle p_1,p_2\rangle = \{m_1 p_1 + m_2 p_2: m_i\in \Z \}$  denote the lattice generated by them.
We say that  $z_1,z_2\in\Z[\ii]$ are \emph{congruent modulo 
    $\langle p_1,p_2\rangle$}, written  $z_1\equiv z_2\pmod{\langle p_1,p_2\rangle}$, if $z_1-z_2 \in \langle p_1,p_2\rangle$.

\begin{definition}\label{def:periodic}
    Let $(a_z)_{z\in\Z[\ii]} $ be a configuration, let $U\subseteq\Z[\ii]$, and let $p_1,p_2\in\Z[\ii]$ be $\Z$-linearly independent.    
    We say that $(a_z)_{z\in\Z[\ii]}$ has \emph{period} $(p_1,p_2)$ in $U$ if for any $z_1,z_2\in U$, whenever $z_1\equiv z_2\pmod{\langle p_1,p_2\rangle}$, we also have $a_{z_1} = a_{z_2}$.

    We say that $(a_z)_{z\in\Z[\ii]}$ has \emph{step-period} $(p_1,p_2)$ in $U$ if whenever $z$ and $z+p_i$  belong to $U$ for some $i$, then $a_z = a_{z+p_i}$.
\end{definition}

\begin{figure}
    \centering
    \includegraphics[width=0.4\linewidth]{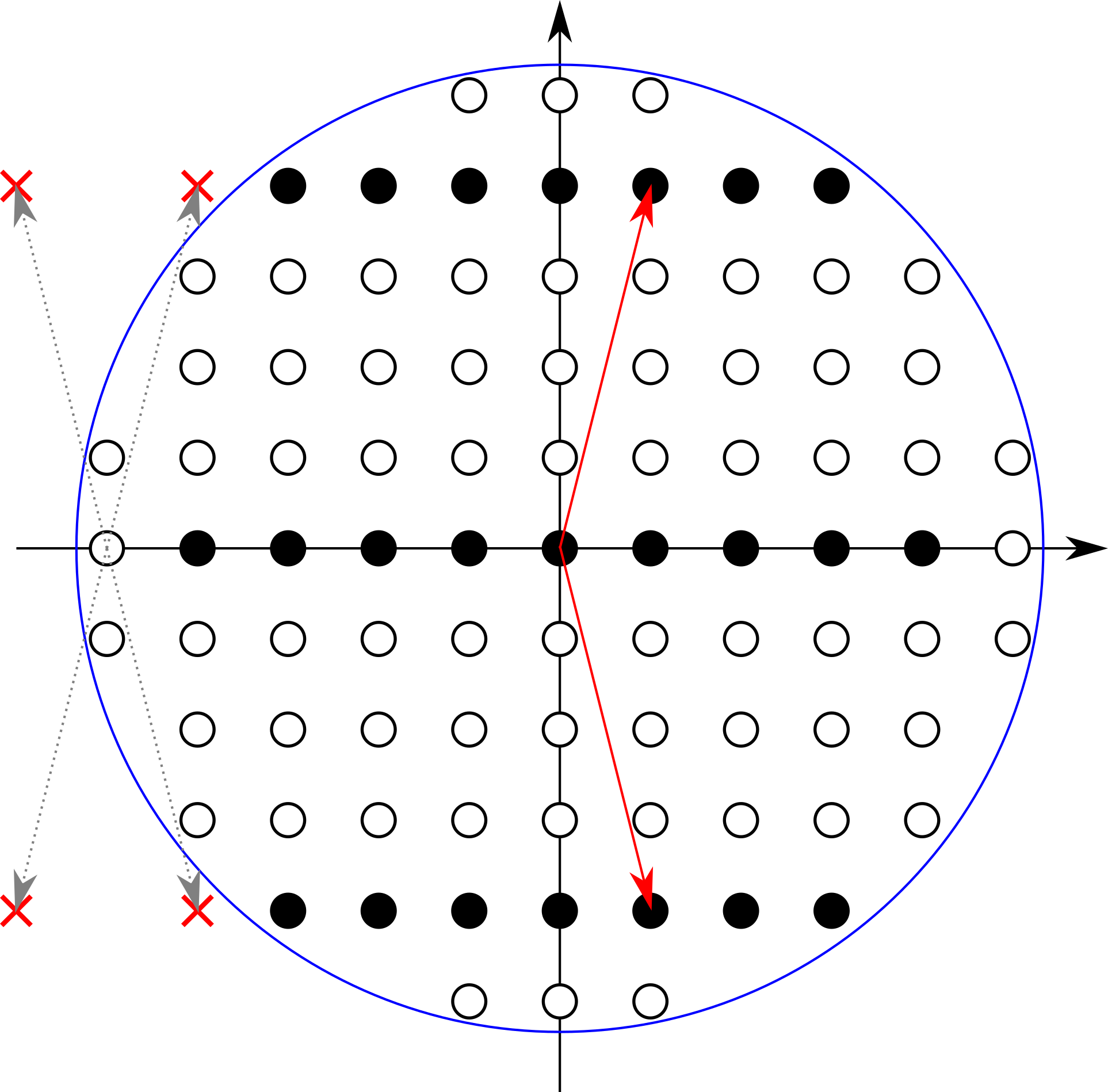}
    \caption{A configuration with step-period $(1+4\ii,1-4\ii)$ on the ball of radius $5 + \epsilon$, for some $\epsilon <1$, around the origin that does not have a period. The values $z=-5$ and $z=1$ are congruent modulo the lattice $\langle 1+4\ii,1-4\ii\rangle$, so the definition of periodicity should require them to be the same color; however, all four neighbours $-5\pm(1\pm 4\ii)$ fall outside the ball, so step-periodicity is not enough to guarantee that this actually happens.}
    \label{fig:step-periodicity-is-not-periodicity}
\end{figure}

It is straightforward to see that if a configuration has period $(p_1,p_2)$ in $U$, then it has step-period $(p_1,p_2)$ in $U$. The converse is true if $U$ is reasonable and we shrink it slightly, as we show in Lemma~\ref{lem:step-lemma}. In particular, for $U=\Z[\ii]$, the two notions are equivalent. For an example of a set which is step-periodic but not periodic, see Figure~\ref{fig:step-periodicity-is-not-periodicity}.

Note that our notion of periodicity is different to that of Hansel and Safer in \cite{Hansel-Safer}. They define $S\subset \Z[\ii]$ to be periodic if there exists $p\in \Z[\ii]$ such that $s\in S$ if and only if $s+zp\in S$ for any $z\in \Z[\ii]$, i.e., that $S$ is a union of cosets of an ideal $(p)$. 
We may even restrict to $p\in \Z$, since if $S$ is a union of cosets mod $(p)$, then it is also a union of
cosets mod $(pq)$, and we may take $q$ to be the complex conjugate of $p$.
Beyond the cosmetic difference that we define periodicity for configurations while they define periodicity for sets, if $U=\Z[\ii] $, one can move between our definition and theirs.  To see this, note that  any lattice $ \langle p_1,p_2\rangle$, where $p_1$ and $p_2$ are  $\Z$-linearly independent, contains elements  $m $ and $n \ii $ with $m$ and $n$  non-zero integers. Setting $p = \lcm(m,n)$,  the lattice  $ \langle p_1,p_2\rangle$ contains any element of the form $ap+bp\ii$, i.e., it contains the ideal $(p)$.  Conversely, being $p$-periodic in the sense of Hansel and Safer  is just
having periods $p$ and $p \ii $ in our first definition.


The proof of the following lemma follows from the triangle inequality.

\begin{lemma} \label{lem:intersection-lemma}
    Let $z_1,z_2\in\C$ and $r_1,r_2\in\R_{>0}$, with $\lvert z_1-z_2\rvert \le r_1 + r_2$. Suppose that neither of the balls $B(z_1,r_1)$ and $B(z_2,r_2)$ is contained in one another. Then their intersection contains a ball of radius $\frac{1}{2}(r_1+r_2-\lvert z_1 - z_2\rvert)$.
\end{lemma}

%

In what follows, we shall study how periodic behaviour extends from one set to another when they have a sufficiently large overlap. As our proofs deal with bounded sets, we have to take into account the fact that the behaviour might not be what is expected near the ``border'' of one such set, and thus we introduce the notation
\[V^{\circ r} \coloneqq \{z\in V \::\: B(z,r)\cap\Z[\ii]\subseteq V \}.\]

The following  two lemmas are reminiscent of  Fine and Wilf's theorem \cite{Fine-Wilf}, which in turn improves Bezout's identity. In both we assume that $p$ and $q$ are multiplicatively independent.

\begin{lemma}[Period transfer lemma]\label{lem:period-transfer-lemma}
    Let $U,V\subseteq\Z[\ii]$ be two sets  such that $(a_z)_{z\in\Z[\ii]}$ has a step-period $(p,q)$ in $U$, and a period $(p',q')$ in $V$. Define the following values:
        \[ \rho_1=\max(\lvert p\rvert,\lvert q\rvert),\quad \rho_2=\max(\lvert p'\rvert,\lvert q'\rvert),\quad R=\rho_1+\rho_2.\]
    If $U\cap V$ contains a ball of radius $R$, then $(a_z)_{z\in\Z[\ii]}$ has a step-period $(p,q)$ in the set $U\cup V^{\circ \rho_1}$.
\end{lemma}

\begin{proof}
    Suppose both $z$ and $z+p$ belong to the set $U\cup V^{\circ \rho_1}$. We are going to show that $a_{z}=a_{z+p}$; this is immediate if both $z$ and $z+p$ are in $U$, so we can assume without loss of generality that either $z$ or $z+p$ belong to $V^{\circ \rho_1}$. Then, both of them belong to $V$; indeed, if we suppose $z\in V^{\circ \rho_1}$, we must have $z+p\in B(z,\rho_1)\cap\Z[\ii]\subseteq V$. The same argument holds for the second case.

    Let $B^*$ be a ball of radius $R$ contained in $U\cap V$, and $B^{**}$ a ball of radius $\rho_2<R$ with the same centre. Note that $B^{**}$ is large enough to contain a complete set of representatives modulo $\langle p',q'\rangle$. Thus, as $z$ is in $V$, there exist $m,n\in\Z$ such that $z+mp'+nq'\in B^{**}\subseteq B^*$. Furthermore, the ball of radius $\rho_1$ centred at $z+mp'+nq'$ is contained in $B^*$, so $z+p+mp'+nq'\in B^*$ as well. Since $B^*\subseteq V$, $a_z=a_{z+mp'+nq'}$ and $a_{z+p}=a_{z+p+mp'+nq'}$ by $(p',q')$-periodicity at $V$. But $B^*$ is also contained in $U$, so $a_{z+mp'+nq'}=a_{(z+mp'+nq')+p}$. We conclude by transitivity that $a_z=a_{z+p}$. The argument for $q$ is similar.
\end{proof}

\begin{lemma}\label{lem:step-lemma}
    Let $U\subseteq\Z[\ii]$ be a ball of radius $R$ such that $(a_z)_{z\in\Z[\ii]}$ has  step-period $(p,q)$ in $U$. Then $(a_z)_{z\in\Z[\ii]}$ has  period $(p,q)$ in $U^{\circ r}$ for $r=|p|+|q|$.
\end{lemma}

\begin{proof}
    We start by defining 
        \begin{equation}\label{eq:rhombus-fundamental-domain}
            \Delta\coloneq\{Ap+Bq\::\: 0\le A,B<1\}\cap\Z[\ii],
        \end{equation}
   which is a complete set of representatives for the lattice $\langle p,q\rangle$, meaning that any element of $\Z[\ii]$ is congruent to some element of $\Delta$.
    
    Let $F\subset\Z[\ii]$ be any connected subset of the Cayley graph of $(\Z[\ii],+)$ with generating set $\{1,\ii\}$, meaning that for any $u,v\in F$ there exists a sequence $u=u_0,\dotsc,u_n=v$ of elements of $F$ satisfying $u_{j+1}-u_j\in\{\pm 1,\pm\ii\}$. If we write $F^{(p,q)}\coloneqq \{ap+bq\::\:a,b\in\Z,a+b\ii\in F\} $, then, for any elements    $z,z'$ in the Minkowski sum $F^{(p,q)}+\Delta$, we have that $z\equiv z'\pmod{\langle p,q\rangle}$ if and only if there exists some sequence $z=z_0,z_1,\dotsc,z_m=z'$ with $z_{j+1}-z_j\in\{\pm p, \pm q\}$ and every $z_j\in F^{(p,q)}+\Delta$, something that follows from the connectedness of $F$. Thus, if $(a_z)_{z\in\Z[\ii]}$ has step-period $(p,q)$ on $F^{(p,q)}+R$, whenever $z\equiv z'\pmod{\langle p,q\rangle}$ we have $a_z=a_{z_1}=a_{z_2}=\cdots=a_{z_m}=a_{z'}$, and hence $(a_z)_{z\in\Z[\ii]}$ has period $(p,q)$ in $F^{(p,q)}+\Delta$ as well, i.e. both properties are equivalent on any set of this form. Naturally, the same holds for any translate of a set of the form $F^{(p,q)}+\Delta$. We have shown that on a set of the form $F^{(p,q)}+\Delta$, step-periodicity is the same as periodicity.

    Now, suppose that $(a_z)_{z\in\Z[\ii]}$ has step-period $(p,q)$ in a ball $U$ of radius $R$. Then, it also has step-period $(p,q)$ in any subset of the form $F^{(p,q)}+\Delta$ contained in $U$, as defined above. Let $F$ be the largest connected subset of $\Z[\ii]$ for which $F^{(p,q)}+ \Delta\subseteq U$; as $r$ is at least the length of either diagonal of $\Delta$, then for any $z\in U^{\circ r}$ there exists some translate  $z'+\Delta$ that contains $z$ with $z'\in F^{(p,q)}$. Thus, $U^{\circ r}\subseteq F^{(p,q)}+\Delta \subseteq U$. As the latter set has step-period $(p,q)$, the middle set has period $(p,q)$, as discussed above. This in turn implies that $U^{\circ r}$ has period $(p,q)$, being a subset of a set with this property. 
\end{proof}



Furthermore, to obtain the result we aim for, we will need every Gaussian integer to be contained in one of the balls $B(x\beta^n,(1-\frac{4}{K})\lvert\beta^n\rvert)$,  with perhaps finitely many exceptions. A sufficient condition for this is that these balls are a cover of $\C$. We state this condition as a lemma:

\begin{lemma}\label{lem:covering_by_balls}
    If $K\ge 28$, then the set $\{B(x\beta^n,(1-\frac{8}{K})\lvert\beta^n\rvert)\::\:x\in\Z[\ii]\}$ is a cover of $\C$ (and thus of $\Z[\ii]$).
\end{lemma}

\begin{proof}
    Geometrically, multiplication by $\beta^n$ corresponds to a rotation followed by a dilation by $\lvert \beta^n\rvert$. Thus, the result is equivalent to the set $\{B(x,1-\frac{8}{K})\::\:x\in\Z[\ii]\}$ being a cover of $\C$. Since the covering radius of the lattice $\Z[\ii]$ is $\frac{\sqrt{2}}{2}$, this will be true for any value of $K$ that ensures that $1-\frac{8}{K}> \frac{1}{\sqrt{2}}$. Hence:
    \[\frac{8}{K}<\frac{2-\sqrt{2}}{2}\implies K>\frac{16}{2-\sqrt{2}}= 8(2+\sqrt{2}).\]
    In particular, $28>8(2+\sqrt{2})$, so these balls cover $\C$ for any value of $K\ge 28$.
\end{proof}

\section{\texorpdfstring{Proof of the main theorem for $\Z[\ii]$}{Proof of the main theorem for $\Z[\ii ]4}}\label{sec:main}

 \begin{lemma}\label{lem:periodic-is-automatic}
  Any eventually periodic configuration $(a_z)_{z\in \Z[\ii]}$ is automatic for any numeration system $(\beta , D)$.  
\end{lemma}

\begin{proof}
We start by recalling some basic facts to ease the exposition. Given any positive integer $p$, every Gaussian integer $z\in\Z[\ii]$ is congruent modulo $p$ to a unique element in the set $\{a+b\ii\::\:0\le a,b<p\}$, which we shall denote by $z\bmod p$; this element may be computed by taking the remainder modulo $p$ of both the real and imaginary parts of $z$, i.e., $a+b\ii\bmod p = (a\bmod p) + (b\bmod p)\ii$. Recall  that for configurations supported on the whole of $\Z[\ii]$ our notion of being periodic in Definition~\ref{def:periodic} is equivalent to the notion of $p$-periodicity introduced by Hansel and Safer, and furthermore we may assume $p$ to be real (more precisely, a positive integer); that is, a configuration $(a_z)_{z\in\Z[\ii]}$ is periodic if and only if there exists $p\in\Z_{>0}$ for which $a_z=a_{z+kp}$ for any $z,k\in\Z[\ii]$. We shall use those two facts to show that a fully periodic configuration is always $(\beta , D)$-automatic for all choices of $D$ and $\beta$. In what follows, we first assume that $(a_z)_{z\in\Z[\ii]}$ is  a periodic configuration, and that $p\in\Z_{>0}$ is chosen such that this configuration is $p$-periodic in the sense of Hansel and Safer.

Note that the sequence $(\beta^j\bmod p)_{j\in\N}$ is eventually periodic,
 so that there exist $m>0$ and $n\geq 0$ such that  $\beta^j\equiv\beta^{j+km}\pmod{p}$ for any $k\ge 0,j\ge n$. In particular, any $\beta^j$ for $j\ge n$ is congruent to some element of $\{\beta^n,\beta^{n+1},\dotsc,\beta^{n+m-1}\}$. We take $n$ and $m$ to be the smallest possible values for these constants.

Given a numeration system $(\beta , D)$, we first define a deterministic finite automaton $\mathcal A$ with state set $S=\{s_{z,j}\::\: z=a+b\ii, 0\le a,b < p, 0\le j < n+m\}$. 
This automaton will be used to compute the remainder of $[w]_\beta \mod p$. The  first index of the state $s_{z,j}$ records the congruence class mod $p$, and the
 second index $j$ tracks which power of $\beta \bmod p$ is concerned, and is directly a function of $|w|$.
The  initial state $s_{0,0}$ and transition function given by:
    \[\delta(s_{z,j},\ell) = \begin{cases}
        s_{z+\beta^j\ell\bmod p,j+1} & \text{if }j<n+m-1, \\
        s_{z+\beta^j\ell\bmod p,n} & \text{if }j=n+m-1,
    \end{cases}\]
    for each $\ell \in D$.
The eventual periodicity of $(\beta^j\bmod{p})$ shows that if $\delta(s_{0,0},w)=s_{z,j}$, then $\beta^{\lvert w\rvert}\equiv\beta^j\pmod{p}$. We claim that from this it follows that $[w]_\beta\equiv z\pmod{p}$, if $w$ is input in reverse reading. The proof is by (structural) induction; the result is immediate when $\lvert w\rvert=1$, as $\delta(s_{0,0},\ell)=s_{\ell\bmod p,1}$. Thus, suppose that the result holds for a given $w\in D^*$, and let $\delta(s_{0,0},\ell w) = \delta(s_{z,j},\ell)=s_{z',j'}$; by the induction hypothesis, $[w]_\beta\equiv z\pmod{p}$. Since $\beta^{\lvert w\rvert}\equiv \beta^j\pmod{p}$, we have:
    \[ [\ell w]_\beta = \beta^{\lvert w\rvert}\ell +[w]_\beta \equiv \beta^j\ell + z \equiv z' \pmod{p}, \]
which proves our claim. This result may be restated as the equality $\delta(s_{0,0},w)=s_{[w]_\beta\bmod{p},j}$, for some $j$ depending only on $\lvert w\rvert$.

As $(a_z)_{z\in\Z[\ii]}$ is $p$-periodic, we have that $z\equiv z'\pmod{p}$ implies that $a_z=a_{z'}$, so in particular $a_z=a_{z\bmod p}$ for any $z\in\Z[\ii]$. If we now define an output function for $\mathcal A$ as $\tau(s_{z,j}) = a_{z\bmod p}$, the corresponding reverse reading DFAO satisfies:
\[\tau(\delta(s_{0,0},w))=\tau(s_{[w]_\beta\bmod p,j})=a_{[w]_\beta \bmod p}=a_{[w]_\beta}, \]
so this is a $(\beta , D)$-automaton that computes $(a_z)_{z\in\Z[\ii]}$.

If $(a_z)_{z\in\Z[\ii]}$ is eventually periodic, there exists a finite set $F\subset\Z[\ii]$ and a fully periodic configuration $(b_z)_{z\in\Z[\ii]}$ such that $a_z=b_z$ for all $z\in\Z[\ii]\setminus F$. Let $\tilde{F}\subset D^*$ be the set of all $w$ for which $[w]_\beta\in F$. If $\tilde{F}$ is finite, then  by   adding finitely many states we can recognize all words of $\tilde{F}$, and define an output  function given by $a_{[w]_\beta}$ whenever $w\in\tilde{F}$, and the original output $b_{[w]_\beta}$ otherwise. If $\tilde{F}$ is infinite (which may happen if $D$ is redundant), then we can use Proposition~\ref{prop:automaticity-equivalent} to convert the numeration system to one where $\tilde{F}$ is finite. The result follows. 

\end{proof}

We can now prove our main theorem. The argument is an  extension of Krebs's argument \cite{Krebs} to the Gaussian numbers.
\begin{main theorem}
Let $\alpha$ and $\beta$ be two multiplicatively independent Gaussian integers with $|\alpha|, |\beta|>1$. Suppose that 
 $\beta$ is not the root of an integer.
Then a configuration is  $\alpha$- and $\beta$-automatic if and only if  it is eventually periodic. 
\end{main theorem}

\begin{remark}\label{rem:best}
The statement of Theorem~\ref{thm:main-gaussian} is the best possible. For,
suppose that $\alpha$ and $\beta$  are both roots of an integer. Then there exist configurations which are both $\alpha$-
 and $\beta$- automatic, but not eventually periodic. To see this we apply \cite[Theorem 4]{Bosma-Fokkink-Krebs}, which tells us that $\mathbb N$ is a
 $\beta$-automatic set if and only if $\beta$ is a root of an integer. The characteristic function of $\mathbb N$ is not eventually periodic. 
   It can be shown that a Semenov-Cobham type theorem exists for such bases.\end{remark}

\begin{proof}[Proof of Theorem~\ref{thm:main-gaussian}] By Lemma~\ref{lem:periodic-is-automatic}, 
    any eventually periodic configuration $(a_z)_{z\in \Z[\ii]}$ is automatic for any base.    
    Conversely, we shall use Proposition~\ref{prop:automaticity-equivalent}, which tells us that an automatic configuration generated by an integral numeration $(\gamma, D)$ is also automatic when generated by a consistent $(\gamma, D_\gamma)$-automaton for a sufficiently enlarged digit set $D_\gamma$ satisfying the conditions of \eqref{eq:expanded_digit_set}. Suppose that $f=(a_z)_{z\in \Z[\ii]}$ is  $\alpha$- and $\beta$-automatic, generated by the finite state machines $\mathcal{A_\alpha}=(S_\alpha,D_\alpha,\delta_\alpha,s_{0,\alpha},A_\alpha,\tau_{\alpha})$ and $\mathcal{A_\beta}=(S_\beta,D_\beta,\delta_\beta,s_{0,\beta},A_{\beta},\tau_{\beta})$ respectively, where $D_\alpha$ and $D_\beta $ are defined as in Equation~\eqref{eq:expanded_digit_set}. 
     Let $\gamma\in \{\alpha, \beta\}$.
     Then by Lemma~\ref{lem:small-numbers-have-short-expansions}, if $z\in B(\gamma^n, |\gamma^n|)\subset  B(0, r|\gamma^n|)$, where $r\geq 2$,  then it will have a short $(\gamma, D_\gamma)$-expansion, specifically  $B(\gamma^n, |\gamma^n|)\subset [D_\gamma^n]_\gamma$. Taking $L_{\gamma s}$ as defined in \eqref{eq:words-ending-at-s}, by Lemma~\ref{lem:shifting-expansions}
        \[a_{x\gamma^n+z} = a_{y\gamma^n+z}\]
    for all $x, y \in [L_{\gamma s}]_\gamma$, all $n \in \N$ and $z \in [D^n_\gamma]_\gamma$.

Recall that  $S_\infty$ is the set of states in $S_\beta$ for which $[L_{\beta, s}]_\beta$ is infinite. If $s\in S_\infty$,
then by Lemma ~\ref{lem:non-collinearity-with-two-automata}
 there exists a state $t=t(s)\in S_\alpha$ for which we may choose three non-collinear Gaussian integers 
 such that each has a $(\alpha, D_\alpha)$-expansion arriving at state $t$ and a $(\beta, D_\beta)$-expansion arriving at state $s$, i.e., 
 $x_{s,t},y_{s,t},z_{s,t}\in [L_{\beta,s}]_\beta\cap [L_{\alpha,t}]_\alpha$. 
 Let
 \[\xi = \max\{\lvert x_{st}\rvert , \lvert y_{st}\rvert, \lvert z_{s,t}\rvert \::\: s \in S_\infty\} + 1,\] 

We  take $\varepsilon=\frac{1}{\xi K}$ in Corollary~\ref{cor:dirichlet}, finding $m,n$ so that
\begin{equation}\label{eq:close-powers-xi-K-version}
    \xi \lvert \alpha^m-\beta^n\rvert<\frac{1}{K}\lvert \beta^n\rvert,
\end{equation} 
where $K$ is a natural number that we will be free to choose when it is relevant. We shall fix these values of $m$ and $n$ from now on. An application of the triangle inequality gives
\begin{equation}\label{eq:alpha-wins-over-beta}
    \lvert\alpha^m\rvert\ge\left(1-\frac{1}{\xi K}\right)\lvert\beta^n\rvert\ge\left(1-\frac{1}{K}\right)\lvert\beta^n\rvert.
\end{equation}

With $x_{s,t},y_{s,t},z_{s,t}$, we may also define two $\R$-linearly independent elements of $\Z[\ii]$ for each $s\in S_\infty$, as follows:
\begin{equation}\label{eq:period-definition}
    p_{s,t}=(y_{s,t}-x_{s,t})(\alpha^m-\beta^n),\quad q_{s,t}=(z_{s,t}-x_{s,t})(\alpha^m-\beta^n).
\end{equation}
We claim that $|p_{s,t}|$ and $|q_{s,t}|$ are small relative to $|\beta^n|$. In particular
\begin{align*}
    \lvert p_{s,t}\rvert &= \lvert y_{s,t} - x_{s,t}\rvert\cdot\lvert\alpha^m-\beta^n\rvert \\
    &\le (\lvert y_{s,t}\rvert + \lvert x_{s,t}\rvert)\cdot \lvert\alpha^m-\beta^n\rvert \\
    &\le 2\xi \lvert\alpha^m-\beta^n\rvert \\
    &\le \frac{2}{K}\lvert\beta^n\rvert;
\end{align*}
similarly for $\lvert q_{s,t}\rvert$.

We want to show that there exist a collection of balls $B$ for which $(a_z)_{z\in\Z[\ii]}$ exhibits a given period $(p_B,q_B)$ in each ball $B$, and such that the balls cover most of $\Z[\ii]$.

\begin{claim}
    $(a_z)_{z\in\Z[\ii]}$ has step-period $(p_{s,t},q_{s,t})$ on the ball $B((x+1)\beta^n,(1-\frac{2}{K})\beta^n)$ for any $x\in [L_{\beta, s}]_\beta$.
\end{claim}



\begin{claimproof}
We only prove the result for $p_{s,t}$; the proof for $q_{s,t}$ is similar. Let $z\in\Z[\ii]$ be any Gaussian integer for which both $z$ and $z+p_{s,t}$ belong to the ball $B(\beta^n,(1-\frac{2}{K})\lvert\beta^n\rvert)$. We start by showing that the number $z-y_{s,t}(\alpha^m-\beta^n)$ is in the ball $B(\alpha^m,\lvert\alpha^m\rvert)\subseteq B(0,r\lvert\alpha^m\rvert)$, and thus has a $(\alpha, D_\alpha)$-expansion of length at most $m$. By an application of the triangle inequality, we obtain:
    \begin{align*}
        \lvert z-y_{s,t}(\alpha^m-\beta^n)-\alpha^m\rvert &\le \lvert z-\beta^n\rvert + \lvert(y_{s,t}+1)(\alpha^m - \beta^n)\rvert \\
        &\le \lvert z-\beta^n\rvert + (\lvert y_{s,t}\rvert + 1) \lvert \alpha^m - \beta^n\rvert \\
        &\le \left(1-\frac{2}{K}\right)\lvert\beta^n\rvert + \xi\cdot\frac{1}{\xi K}\lvert\beta^n\rvert \\
        & \le \left(1-\frac{1}{K}\right)\lvert\beta^n\rvert \\
        & \le \lvert\alpha^m\rvert,
    \end{align*}
    where the last inequality follows from Equation~\eqref{eq:alpha-wins-over-beta}.

    Now, since $x\in[L_{\beta,s}]_\beta$, it must have at least one $(\beta , D_\beta)$-expansion that ends at state $s$, so we can apply Lemma~\ref{lem:shifting-expansions}, to get the following:
    \begin{align*}
        a_{x\beta^n+z} &=a_{y_{s,t}\beta^n+z} & &\text{since } x,y_{s,t}\in[L_{\beta,s}]_\beta \text{ and}\\ & & & \text{$z$ has a $(\beta, D_\beta)$-expansion of length${}\le n$,} \\
        & =a_{y_{s,t}\alpha^m+z-y_{s,t}(\alpha^m-\beta^n)} & & \text{by adding and substracting $y_{s,t}\alpha^m$,}\\
        &=a_{x_{s,t}\alpha^m+z-y_{s,t}(\alpha^m-\beta^n)} & & \text{since }x_{s,t},y_{s,t}\in [L_{\alpha,t}]_\alpha \text{ and}\\ & & & \text{$z-y_{s,t}(\alpha^m-\beta^n)$ has a $(\alpha, D_\alpha)$-expansion of length${}\le m$,} \\
        &=a_{x_{s,t}\beta^n + z +p_{s,t}} & &\text{by adding and substracting $x_{s,t}\beta^n$ and} \\ & & &  \text{replacing $(x_{s,t}-y_{s,t})(\alpha^m-\beta^n )$ by $p_{s,t}$,} \\
        &=a_{x\beta^n+z+p_{s,t}} & &\text{since $x_{s,t},x\in [L_{\beta,s}]_\beta$ and} \\ & & & \text{$z+p_{s,t}$ has a $(\beta, D_\beta)$-expansion of length${}\le n$},
    \end{align*}
    and this shows, thus, that $a_{u}=a_{u+p_{s,t}}$ for $u=x\beta^n+z\in B((x+1)\beta^n,(1-\frac{2}{K})\lvert\beta^n\rvert)$. This ends the proof of the claim.
\end{claimproof}

Therefore, by Lemma~\ref{lem:step-lemma}, $(a_z)_{z\in\Z[\ii]}$ has period $(p_{s,t},q_{s,t})$ on the ball $B((x+1)\beta^n,(1-\frac{6}{K})\beta^n)$ for any $x\in [L_{\beta, s}]_\beta$.

All but finitely many Gaussian integers have a $(\beta, D_\beta)$-expansion that ends at a state in $S_\infty$. For each Gaussian integer $u$ not in this finite set $F$, we fix $(p_{u},q_u)$ as above so that $(a_z)_{z\in\Z[\ii]}$ has period $(p_{u},q_{u})$ on the ball $B_u\coloneqq B((u+1)\beta^n,(1-
\frac{6}{K})\lvert\beta^n\rvert)$. We also introduce the notation $C_u\coloneqq B_{u}^{\circ (2 \lvert\beta^n\rvert/K)}$ to denote a ball of radius $(1- 
\frac{8}{K}
)\lvert\beta^n\rvert$ with the same centre.

For any two balls $C_u$ and $B_{u+z}$ with $z\in\{\pm 1,\pm\ii\}$, the distance between their centres will be $\lvert\beta^n\rvert$. Thus, Lemma~\ref{lem:intersection-lemma}, with $r_1= (1- 
\frac{6}{K})\beta^n$ and $r_2= (1-  \frac{8}{K})\beta^n$, guarantees that their intersection contains a ball of radius 
$\frac{1}{2}(2-\frac{14}{K}-1)\lvert \beta^n\rvert=\frac{K-14}{2K}\lvert \beta^n\rvert$
. 

As  $|p_{u}|, |q_{u}|, |p_{u+z}|, |q_{u+z}|$ are each bounded by $\frac{2}{K}\lvert\beta^n\rvert$,  
then the intersection $C_u \cap B_{u+z}$ contains a ball of radius as required by Lemma~\ref{lem:period-transfer-lemma}, if $K$ is large enough. If neither $u$ nor $u+z$ belong to $F$ we can apply Lemma~\ref{lem:period-transfer-lemma}, taking $U=C_n$ and $V=B_{n+z}$, if the intersection of these balls contains a ball of radius at least $\frac{4}{K}\lvert\beta^n\rvert$. Hence, we have the inequality:
    \[
     \frac{K-14}{2K}\lvert\beta^n\rvert 
    \ge \frac{4}{K}\lvert\beta^n\rvert\implies K\ge  22,\]
so any sufficiently large choice of $K$ allows us to apply this result and to conclude that $(a_z)$ has a step-period $(p_u,q_u)$ in $C_u\cup B_{u+z}^{\circ(2\lvert\beta^n\rvert/K)}=C_u\cup C_{u+z}$.

Lemma~\ref{lem:covering_by_balls} ensures that, whenever $K\geq 28$, the balls $\{C_u\::\: u \not\in F\} $ cover the whole of  $\Z[\ii]$ except for some subset of $\left(\bigcup_{u\in F} C_u\right)\cap\Z[\ii]$, which is a finite set, as it is a bounded subset of $\Z[\ii]$. By enlarging the set $F$ if necessary, we may ensure that the induced subgraph of the Cayley graph of $(\Z[\ii],+)$ with generating set $\{\pm1,\pm\ii\}$ is connected.

Fix any $u^*\in \Z[\ii]\setminus F$, and let $v\in\Z[\ii]\setminus F$ be any other element of this set. Let $u^*=u_0,u_1,\dotsc,u_n=v$ be a path between $u^*$ and $v$; since we see a step-period $(p_{u^*},q_{u^*})$ in $C_{u^*}$ and a period $(p_{u_1},q_{u_1})$ in $B_{u_1}$, an application of Lemma~\ref{lem:period-transfer-lemma} shows that we observe a step-period $(p_{u^*},q_{u^*})$ in $C_{u_1}$; proceeding inductively, we show that a step-period $(p_{u^*},q_{u^*})$ in $C_{u_j}$ induces the same step-period on $C_{u_{j+1}}$. This applies, in particular, to $C_v$. As this is true for any $v\in\Z[\ii]\setminus F$, we conclude that $(a_z)_{z\in\Z[\ii]}$ has a step-period $(p_{u^*},q_{u^*})$ in $\left(\bigcup_{u\in\Z[\ii]\setminus F} C_u\right)\cap\Z[\ii]$, which, as discussed before, equals the whole of $\Z[\ii]$ bar a finite set.
The result follows.
%
%
%
\end{proof}

As a corollary, we answer a question of \cite{Bosma-Fokkink-Krebs}. 
\begin{corollary}
Suppose that $U\subset \Z[\ii]$ is a $\beta$-automatic set, where $|\beta|>1$. If $\overline U = U$, and $U$ is not eventually periodic, then $\beta$ is a root of an integer.
\end{corollary}
\begin{proof} Suppose that $U$ is $(\beta, D)$-automatic. If needed,  using Proposition~\ref{prop:automaticity-equivalent}, we can enlarge $D$ so that it is closed under conjugation. Since $U=\overline U$, we conclude that $U$ is $(\bar \beta, D)$-automatic.  Since $U$ is not eventually periodic,  $\beta$ and $\bar\beta$ are multiplicatively dependent. Thus there exist $m, n$ such that  $\beta^n=\bar \beta^m$ and so
$|\beta|^n=|\bar \beta|^m$. This implies that  $m=n$, since $|\beta|>1$. Thus 
$\beta^n=\bar \beta^n$, so that
$\beta^n\in\mathbb R$.
\end{proof}

Hansel and Safer were only interested in ``natural'' integral Gaussian numerations $(\beta,D)$, i.e., those where  $D\subset \mathbb N$. 
The fact   that  $D$  is an integral digit set implies that $\beta\in \{-n\pm i : n\geq 1 \}$ by Part (4) of Theorem~\ref{thm:numerations-exist}. Thus given two multiplicatively independent  $\alpha, \beta$ with integral digit sets, at least one of them is not a root of an integer. 
   The following result now follows immediately from Theorem~\ref{thm:main-gaussian}.

\begin{corollary}
Let $\alpha$ and $\beta$ be two multiplicatively independent Gaussian integers with $|\alpha|, |\beta|>1$. Suppose that there exist digit sets
 $D_\alpha$ and $D_\beta$ consisting of natural numbers such that $(\alpha,D_\alpha)$ and $(\beta,D_\beta)$ are integral numeration systems.
Then a configuration is  $\alpha$- and $\beta$-automatic  if and only if  it is eventually periodic. 
\end{corollary}

Finally, we mention the work of Cabezas and Petite \cite{Cabezas-Petite} on constant shape substitutions. Via the representation of $\Z[\ii]$ with the ring of integer matrices of the form $\begin{psmallmatrix} a & -b\\ b& a \end{psmallmatrix}$, a $\beta$-automatic configuration is also the coding of a fixed point of a constant shape substitution, determined by an integer $2\times 2$ matrix $M$ and  fundamental domain $D$ for $M(\Z^2)$. In the case where $M$ represents multiplication by $\beta\in \Z[\ii]$, the configurations that Cabezas and Petite generate are $(\beta,D)$-automatic configurations. Our work therefore applies to their configurations, and in particular it means that  if $\beta=a+ib$ is not the root of an integer, then their nontrivial $(\begin{psmallmatrix} a & -b\\ b& a \end{psmallmatrix},D)$-configurations  cannot be replicated by $(k,k)$-automatic configurations, for $k\in \mathbb N$.

In conclusion, it is natural to wonder what happens for other rings of integers. Natural candidates include $\Z[\zeta]$, where $\zeta$ is a root of unity,  quadratic integer rings, or Euclidean domains.

\bibliographystyle{plain}
\bibliography{references}

@article {Coven-Dykstra-Lemasurier,
    AUTHOR ={Coven, Ethan M. and Dykstra, Andrew and Lemasurier, Michelle},
     TITLE = {A short proof of a theorem of {C}obham on substitutions},
   JOURNAL = {Rocky Mountain J. Math.},
  FJOURNAL = {The Rocky Mountain Journal of Mathematics},
    VOLUME = {44},
      YEAR = {2014},
    NUMBER = {1},
     PAGES = {19--22},
      ISSN = {0035-7596},
   MRCLASS = {68Q70 (37B10)},
  MRNUMBER = {3216006},
MRREVIEWER = {Michael H. Schraudner},
       DOI = {10.1216/RMJ-2014-44-1-19},
       URL = {https://doi-org.ezproxy.library.qmul.ac.uk/10.1216/RMJ-2014-44-1-19},
}

@inproceedings {HS-2022,
    AUTHOR = {Hieronymi, Philipp and Schulz, Christian},
     TITLE = {A strong version of {C}obham's theorem},
 BOOKTITLE = {S{TOC} '22---{P}roceedings of the 54th {A}nnual {ACM} {SIGACT}
              {S}ymposium on {T}heory of {C}omputing},
     PAGES = {1172--1179},
 PUBLISHER = {ACM, New York},
      YEAR = {[2022] \copyright 2022},
      ISBN = {978-1-4503-9264-8},
   MRCLASS = {03F30},
  MRNUMBER = {4490070},
       DOI = {10.1145/3519935.3519958},
       URL = {https://doi-org.ezproxy.library.qmul.ac.uk/10.1145/3519935.3519958},
}

@article {Hosseini-Yassawi,
    AUTHOR = {Hosseini, Maryam and Yassawi, Reem},
     TITLE = {Obstacles to topological factoring of {T}oeplitz shifts},
   JOURNAL = {Discrete Contin. Dyn. Syst.},
  FJOURNAL = {Discrete and Continuous Dynamical Systems},
    VOLUME = {46},
      YEAR = {2026},
     PAGES = {413--432},
      ISSN = {1078-0947},
   MRCLASS = {37B10 (37B05)},
  MRNUMBER = {4973347},
       DOI = {10.3934/dcds.2025105},
       URL = {https://doi-org.ezproxy.library.qmul.ac.uk/10.3934/dcds.2025105},
}

@article {Adam-Fav-1,
    AUTHOR = {Adamczewski, Boris and Faverjon, Colin},
     TITLE = {Mahler's method in several variables and finite automata},
   JOURNAL = {Annals of Mathematics},
      VOLUME = {to appear},
       URL = {https://arxiv.org/abs/2012.08283},
}

@article {Adam-Fav-2,
    AUTHOR = {Adamczewski, Boris and Faverjon, Colin},
     TITLE = {Addendum to: Mahler's method in several variables and finite automata},
   JOURNAL = {Annals of Mathematics},
      VOLUME = {to appear},
       URL = {https://arxiv.org/abs/2407.18578},
}

@incollection {BHMV1,
    AUTHOR = {Bruy\`ere, V\'{e}ronique and Hansel, Georges and Michaux, Christian
              and Villemaire, Roger},
     TITLE = {Logic and {$p$}-recognizable sets of integers, (and erratum)},
      NOTE = {Journ\'{e}es Montoises (Mons, 1992)},
   JOURNAL = {Bull. Belg. Math. Soc. Simon Stevin},
  FJOURNAL = {Bulletin of the Belgian Mathematical Society. Simon Stevin},
    VOLUME = {1},
      YEAR = {1994},
    NUMBER = {2 and 4},
     PAGES = {191--238, and 577},
      ISSN = {1370-1444},
   }

@article {Fine-Wilf,
    AUTHOR = {Fine, N. J. and Wilf, H. S.},
     TITLE = {Uniqueness theorems for periodic functions},
   JOURNAL = {Proc. Amer. Math. Soc.},
  FJOURNAL = {Proceedings of the American Mathematical Society},
    VOLUME = {16},
      YEAR = {1965},
     PAGES = {109--114},
      ISSN = {0002-9939},
   MRCLASS = {42.99},
  MRNUMBER = {174934},
MRREVIEWER = {A. Devinatz},
       DOI = {10.2307/2034009},
       URL = {https://doi-org.ezproxy.library.qmul.ac.uk/10.2307/2034009},
}

@book {Knuth:1988,
    AUTHOR = {Knuth, Donald E.},
     TITLE = {The art of computer programming. {V}ol. 2},
      NOTE = {Seminumerical algorithms,
              Third edition [of MR0286318]},
 PUBLISHER = {Addison-Wesley, Reading, MA},
      YEAR = {1998},
     PAGES = {xiv+762},
      ISBN = {0-201-89684-2},
   MRCLASS = {68-02 (65-XX)},
  MRNUMBER = {3077153},
}

@article{Penney,
 AUTHOR = {Penney, Walter},
     TITLE = {A 'binary' system for complex numbers},
year = {1965},
issue_date = {April 1965},
publisher = {Association for Computing Machinery},
volume = {12},
number = {2},
issn = {0004-5411},
journal = {J. ACM}
}

@article {Cabezas-Leroy,
    AUTHOR = {Cabezas, Christopher and Leroy, Julien},
     TITLE = {Decidability of the isomorphism problem between
              multidimensional substitutive subshifts},
   JOURNAL = {Ergodic Theory Dynam. Systems},
  FJOURNAL = {Ergodic Theory and Dynamical Systems},
    VOLUME = {45},
      YEAR = {2025},
    NUMBER = {7},
     PAGES = {2054--2094},
      ISSN = {0143-3857},
   MRCLASS = {37B10 (37A35 37B51 68R15)},
  MRNUMBER = {4915539},
       DOI = {10.1017/etds.2024.87},
       URL = {https://doi-org.ezproxy.library.qmul.ac.uk/10.1017/etds.2024.87},
}

@article {Adamczewski-Bell-2011,
    AUTHOR = {Adamczewski, Boris and Bell, Jason},
     TITLE = {An analogue of {C}obham's theorem for fractals},
   JOURNAL = {Trans. Amer. Math. Soc.},
  FJOURNAL = {Transactions of the American Mathematical Society},
    VOLUME = {363},
      YEAR = {2011},
    NUMBER = {8},
     PAGES = {4421--4442},
      ISSN = {0002-9947},
   MRCLASS = {28A80 (11B85)},
  MRNUMBER = {2792994},
MRREVIEWER = {Tomas Persson},
       DOI = {10.1090/S0002-9947-2011-05357-2},
       URL = {https://doi-org.ezproxy.library.qmul.ac.uk/10.1090/S0002-9947-2011-05357-2},
}

@article {Bes-2000,
    AUTHOR = {B\`es, Alexis},
     TITLE = {An extension of the {C}obham-{S}em\"{e}nov theorem},
   JOURNAL = {J. Symbolic Logic},
  FJOURNAL = {The Journal of Symbolic Logic},
    VOLUME = {65},
      YEAR = {2000},
    NUMBER = {1},
     PAGES = {201--211},
      ISSN = {0022-4812},
   MRCLASS = {03D05 (11R06 68Q45)},
  MRNUMBER = {1782115},
MRREVIEWER = {V\'{e}ronique Bruy\`ere},
       DOI = {10.2307/2586532},
       URL = {https://doi-org.ezproxy.library.qmul.ac.uk/10.2307/2586532},
}

@article {Adamczewski-Bell-2008,
    AUTHOR = {Adamczewski, Boris and Bell, Jason},
     TITLE = {Function fields in positive characteristic: expansions and
              {C}obham's theorem},
   JOURNAL = {J. Algebra},
  FJOURNAL = {Journal of Algebra},
    VOLUME = {319},
      YEAR = {2008},
    NUMBER = {6},
     PAGES = {2337--2350},
      ISSN = {0021-8693},
   MRCLASS = {11B85 (11R58 13F25)},
  MRNUMBER = {2388308},
MRREVIEWER = {David Adam},
       DOI = {10.1016/j.jalgebra.2007.06.039},
       URL = {https://doi-org.ezproxy.library.qmul.ac.uk/10.1016/j.jalgebra.2007.06.039},
}

@article {Durand-2011,
    AUTHOR = {Durand, Fabien},
     TITLE = {Cobham's theorem for substitutions},
   JOURNAL = {J. Eur. Math. Soc. (JEMS)},
  FJOURNAL = {Journal of the European Mathematical Society (JEMS)},
    VOLUME = {13},
      YEAR = {2011},
    NUMBER = {6},
     PAGES = {1799--1814},
      ISSN = {1435-9855},
   MRCLASS = {68Q45 (11A67)},
  MRNUMBER = {2835330},
MRREVIEWER = {Michel Rigo},
       DOI = {10.4171/JEMS/294},
       URL = {https://doi-org.ezproxy.library.qmul.ac.uk/10.4171/JEMS/294},
}

@article {Bell-2007,
    AUTHOR = {Bell, Jason P.},
     TITLE = {A generalization of {C}obham's theorem for regular sequences},
   JOURNAL = {S\'{e}m. Lothar. Combin.},
  FJOURNAL = {S\'{e}minaire Lotharingien de Combinatoire},
    VOLUME = {54A},
      YEAR = {2005/07},
     PAGES = {Art. B54Ap. 15},
   MRCLASS = {11B85},
  MRNUMBER = {2223028},
MRREVIEWER = {Val\'{e}rie Berth\'{e}},
}

@article{vdW,
	author = {van der Waerden, B. L.},
		journal = {  Nieuw Arch. Wisk.    },
		number = {2},
	pages = {212--216},
	title = {  Beweis einer {B}audet'schen Vermutung},
		volume = {15},
	year = {1928},
	}

@article{SS-2019,
	author = {Sch\"{a}fke, Reinhard and Singer, Michael},
	doi = {10.4171/JEMS/891},
	fjournal = {Journal of the European Mathematical Society (JEMS)},
	issn = {1435-9855},
	journal = {J. Eur. Math. Soc. (JEMS)},
	mrclass = {39A05 (34A30 34M03 39A13 39A45)},
	mrnumber = {3985611},
	mrreviewer = {P. W. Eloe},
	number = {9},
	pages = {2751--2792},
	title = {Consistent systems of linear differential and difference equations},
	url = {https://doi-org.ezproxy.library.qmul.ac.uk/10.4171/JEMS/891},
	volume = {21},
	year = {2019}}

@article {Muchnik,
    AUTHOR = {Muchnik, I.},
     TITLE = {Definable criterion for definability in Presburger Arithmetic and
its application},
   JOURNAL = { Institute of New Technologies (preprint in Russian)},
        YEAR = {1991},}

@article {Katai-Szabo,
    AUTHOR = {K\'{a}tai, I. and Szab\'{o}, J.},
     TITLE = {Canonical number systems for complex integers},
   JOURNAL = {Acta Sci. Math. (Szeged)},
  FJOURNAL = {Acta Universitatis Szegediensis. Acta Scientiarum
              Mathematicarum},
    VOLUME = {37},
      YEAR = {1975},
    NUMBER = {3-4},
     PAGES = {255--260},
      ISSN = {0001-6969},
   MRCLASS = {10B35 (10A10)},
  MRNUMBER = {389759},
MRREVIEWER = {J. M. Gandhi},
}

@article {Bosma-Fokkink-Krebs,
    AUTHOR = {Bosma, Wieb and Fokkink, Robbert and Krebs, Thijmen},
     TITLE = {On automatic subsets of the {G}aussian integers},
   JOURNAL = {Indag. Math. (N.S.)},
  FJOURNAL = {Koninklijke Nederlandse Akademie van Wetenschappen.
              Indagationes Mathematicae. New Series},
    VOLUME = {28},
      YEAR = {2017},
    NUMBER = {1},
     PAGES = {32--37},
      ISSN = {0019-3577},
   MRCLASS = {11B85 (68Q45)},
  MRNUMBER = {3597031},
MRREVIEWER = {Eric Rowland},
       DOI = {10.1016/j.indag.2016.11.003},
       URL = {https://doi-org.ezproxy.library.qmul.ac.uk/10.1016/j.indag.2016.11.003},
}

@article {Krebs,
    AUTHOR = {Krebs, Thijmen J. P.},
     TITLE = {A more reasonable proof of {C}obham's theorem},
   JOURNAL = {Internat. J. Found. Comput. Sci.},
  FJOURNAL = {International Journal of Foundations of Computer Science},
    VOLUME = {32},
      YEAR = {2021},
    NUMBER = {2},
     PAGES = {203--207},
      ISSN = {0129-0541},
   MRCLASS = {68Q45 (11B85)},
  MRNUMBER = {4218824},
       DOI = {10.1142/S0129054121500118},
       URL = {https://doi-org.ezproxy.library.qmul.ac.uk/10.1142/S0129054121500118},
}

@book {Lothaire,
    AUTHOR = {Lothaire, M.},
     TITLE = {Algebraic combinatorics on words},
    SERIES = {Encyclopedia of Mathematics and its Applications},
    VOLUME = {90},
      NOTE = {A collective work by Jean Berstel, Dominique Perrin, Patrice
              Seebold, Julien Cassaigne, Aldo De Luca, Steffano Varricchio,
              Alain Lascoux, Bernard Leclerc, Jean-Yves Thibon, Veronique
              Bruyere, Christiane Frougny, Filippo Mignosi, Antonio Restivo,
              Christophe Reutenauer, Dominique Foata, Guo-Niu Han, Jacques
              Desarmenien, Volker Diekert, Tero Harju, Juhani Karhumaki and
              Wojciech Plandowski,
              With a preface by Berstel and Perrin},
 PUBLISHER = {Cambridge University Press, Cambridge},
      YEAR = {2002},
     PAGES = {xiv+504},
      ISBN = {0-521-81220-8},
   MRCLASS = {68R15 (05A99 05E10 20M35)},
  MRNUMBER = {1905123},
MRREVIEWER = {Patrice S\'{e}\'{e}bold},
       DOI = {10.1017/CBO9781107326019},
       URL = {https://doi-org.ezproxy.library.qmul.ac.uk/10.1017/CBO9781107326019},
}

@article {Cobham-1969,
    AUTHOR = {Cobham, Alan},
     TITLE = {On the base-dependence of sets of numbers recognizable by
              finite automata},
   JOURNAL = {Math. Systems Theory},
  FJOURNAL = {Mathematical Systems Theory. An International Journal on
              Mathematical Computing Theory},
    VOLUME = {3},
      YEAR = {1969},
     PAGES = {186--192},
      ISSN = {0025-5661},
   MRCLASS = {94.40 (02.00)},
  MRNUMBER = {250789},
MRREVIEWER = {J. Hartmanis},
       DOI = {10.1007/BF01746527},
       URL = {https://doi-org.ezproxy.library.qmul.ac.uk/10.1007/BF01746527},
}

@article {Hansel-Safer,
    AUTHOR = {Hansel, Georges and Safer, Taoufik},
     TITLE = {Vers un th\'{e}or\`eme de {C}obham pour les entiers de {G}auss},
   JOURNAL = {Bull. Belg. Math. Soc. Simon Stevin},
  FJOURNAL = {Bulletin of the Belgian Mathematical Society. Simon Stevin},
    VOLUME = {10},
      YEAR = {2003},
    NUMBER = {suppl.},
     PAGES = {723--735},
      ISSN = {1370-1444},
   MRCLASS = {68R15 (11B85 68Q45)},
  MRNUMBER = {2073023},
       URL = {http://projecteuclid.org.ezproxy.library.qmul.ac.uk/euclid.bbms/1074791328},
}

@article{Cabezas-Petite,
title = {Large normalizers of $  \mathbb{Z}^{d} $-odometer systems and realization on substitutive subshifts},
journal = {Discrete and Continuous Dynamical Systems},
volume = {44},
number = {12},
pages = {3848-3877},
year = {2024},
issn = {1078-0947},
doi = {10.3934/dcds.2024080},
url = {https://www.aimsciences.org/article/id/667142d65a42b314c5bff304},
author = {Christopher Cabezas and Samuel Petite}
}

@article {Semenov,
    AUTHOR = {Semenov, A. L.},
     TITLE = {The {P}resburger nature of predicates that are regular in two
              number systems},
   JOURNAL = {Sibirsk. Mat. \v{Z}.},
  FJOURNAL = {Akademija Nauk SSSR. Sibirskoe Otdelenie. Sibirski\u{\i}
              Matemati\v{c}eski\u{\i} \v{Z}urnal},
    VOLUME = {18},
      YEAR = {1977},
    NUMBER = {2},
     PAGES = {403--418, 479},
      ISSN = {0037-4474},
   MRCLASS = {02F43},
  MRNUMBER = {450050},
MRREVIEWER = {Constantin Milici},
}

@article {Durand-2008,
    AUTHOR = {Durand, F.},
     TITLE = {Cobham-{S}emenov theorem and {$\Bbb N^d$}-subshifts},
   JOURNAL = {Theoret. Comput. Sci.},
  FJOURNAL = {Theoretical Computer Science},
    VOLUME = {391},
      YEAR = {2008},
    NUMBER = {1-2},
     PAGES = {20--38},
      ISSN = {0304-3975},
   MRCLASS = {68Q45},
  MRNUMBER = {2381349},
MRREVIEWER = {Michel Rigo},
       DOI = {10.1016/j.tcs.2007.10.027},
       URL = {https://doi-org.ezproxy.library.qmul.ac.uk/10.1016/j.tcs.2007.10.027},
}

@article{mihailescu2004,
  title={Primary cyclotomic units and a proof of {C}atalan's conjecture},
  author={Mih{\u{a}}ilescu, Preda},
  journal={J. Reine Angew. Math.},
  volume={572},
  pages={167–195},
  year={2004},
  publisher={Walter de Gruyter GmbH \& Co. KG Berlin, Germany}
}

@article{DDG,
  title={Complex arithmetic},
  author={Davio, M and Deschamps, JP and Gossart, C},
  journal={Philips MBLE Research Lab. Report},
  volume={369},
  year={1978}
}

@book{niven,
  title={Irrational numbers},
  author={Niven, Ivan},
  series={Carus Mathematical Monographs},
  number={11},
  year={2005},
  publisher={Cambridge University Press}
}

\end{document}